\documentclass[article,12pt]{amsart}

\usepackage[T1]{fontenc}
\usepackage{inputenc}
\usepackage[francais,english]{babel}
\usepackage{amsmath}
\usepackage{latexsym}
\usepackage{marvosym}
\usepackage{amsfonts}
\usepackage{amsthm}
\usepackage{float}
\usepackage{color}
\usepackage{epsfig}
\usepackage{graphicx}

\newcommand{\dR}{\mathbb{R}}
\newcommand{\dE}{\mathbb{E}}
\newcommand{\dP}{\mathbb{P}}

\newcommand{\cG}{\mathcal{G}}
\newcommand{\cM}{\mathcal{M}}
\newcommand{\cN}{\mathcal{N}}

\newcommand{\cP}{\mathcal{P}}
\newcommand{\cQ}{\mathcal{Q}}
\newcommand{\cR}{\mathcal{R}}
\newcommand{\cS}{\mathcal{S}}

\newcommand{\cF}{\mathcal{F}}

\newcommand{\veps}{\varepsilon}
\newcommand{\rI}{\mathrm{I}}
\newcommand{\wh}{\widehat}

\newcommand{\ind}{\mbox{1}\kern-.25em \mbox{I}}
\font\calcal=cmsy10 scaled\magstep1
\def\build#1_#2^#3{\mathrel{\mathop{\kern 0pt#1}\limits_{#2}^{#3}}}
\def\liml{\build{\longrightarrow}_{}^{{\mbox{\calcal L}}}}

\def\videbox{\mathbin{\vbox{\hrule\hbox{\vrule height1ex \kern.5em
\vrule height1ex}\hrule}}}

\numberwithin{equation}{section}
\theoremstyle{plain}
\newtheorem{thm}{Theorem}[section]
\newtheorem{rem}{Remark}[section]

\email{philippe.fraysse@math.u-bordeaux1.fr}

\begin{document}
\title[Recursive estimation in a class of models of deformation]
{Recursive estimation in a class of models of deformation}
\author{Philippe Fraysse}
\dedicatory{\normalsize Universit\'e de Bordeaux}
\address{Universit\'e Bordeaux 1, Institut de Math\'ematiques de Bordeaux, UMR CNRS 5251, and INRIA Bordeaux, team ALEA,
	351 cours de la lib\'eration, 33405 Talence cedex, France.}

\begin{abstract}
The paper deals with the statistical analysis of several data sets associated with shape invariant models with different translation, height and scaling parameters. We propose to estimate these parameters together with the common shape function. Our approach extends the recent work of Bercu and Fraysse to multivariate shape invariant models. We propose a very efficient Robbins-Monro procedure for the estimation of the translation parameters and we use these estimates in order to evaluate scale parameters. The main pattern is estimated by a weighted Nadaraya-Watson recursive estimator. We provide almost sure convergence and asymptotic normality for all estimators. Finally, we illustrate the convergence of our estimation procedure on simulated data as well as on real ECG data.  
\end{abstract}

\maketitle


\section{INTRODUCTION}


Statistic analysis of models with periodic data is a mathematical field of great interest. Indeed, a detailed analysis of such models enables us to have a satisfactory approximation of real life phenomena. In particular, 
SEMOR models \cite{Kneip88} are often used to describe a large number of phenomena as Meteorology \cite{Vimond10}, road traffic \cite{CastilloLoubes09}, \cite{GamboaLoubes07} or children growth \cite{Gasser95}. Here, we choose to focus our attention on a particular class of these models called shape invariant models, introduced by Lawton \textit{et al.} \cite{Lawton}. Periodic shape invariant models are semiparametric regression models with
an unknown periodic shape function. In this paper, we consider several data sets associated with a common shape function and differing from each other by three parameters, a translation, a height and a scale. Formally, we are interested in the following shape invariant model
\begin{equation}
\label{Sempar}
Y_{i,j}=a_{j}f(X_{i}-\theta_{j})+v_{j}+\varepsilon_{i,j},
\end{equation}
where
$1\leq{j}\leq{p}$ and $1\leq{i}\leq{n}$, the common shape function $f$ is periodic and the variables $X_{i}$ are random, independent and of the same law. The classical approaches to estimate the different parameters of the model are to minimize the least-squares or to maximize the likelihood of the model when the law of $\varepsilon_{i,j}$ is known. Here, we propose a new recursive estimation procedure which requires very few assumptions and is really easy to implement.
The case where the $X_{i}$ are equi-distributed deterministic variables has been considered in \cite{GamboaLoubes07}, \cite{Trigano11} or \cite{Vimond10}. When $p=1$, $a_{1}=1$ and $v_{1}=0$, Bercu and Fraysse \cite{BF10} propose a recursive method to estimate the translation parameter $\theta_{1}$. 
In this paper, we significantly extend their results as we are now able to estimate, whatever the value of the dimension parameter $p$ is, the height parameter $v$, the translation parameter $\theta$ and the scale parameter $a$, respectively given by
 \begin{equation}
 v=\begin{pmatrix}v_{1}\\
\vdots\\
v_{p}\end{pmatrix},
\hspace{1cm}
 \theta=\begin{pmatrix}\theta_{1}\\
\vdots\\
\theta_{p}\end{pmatrix},
\hspace{1cm}
 a=\begin{pmatrix}a_{1}\\
\vdots\\
a_{p}\end{pmatrix}.
 \end{equation}
 \par
 Our first goal is to estimate the translation parameter $\theta$. Estimation of shifts has lots of similarities with curve registration and alignment problems \cite{RamsaySilver05}. Analysis of ECG curves \cite{Trigano11} or the study of traffic data \cite{CastilloLoubes09}, \cite{GamboaLoubes07} fall into this framework. In \cite{Gasser95}, Gasser and Kneip propose to estimate the shifts by aligning the maxima of the curves, their position being estimated by the zeros  of the derivative of a kernel estimate. In the case where $a_{j}=1$ and $v_{j}=0$, Gamboa \textit{et al.} \cite{GamboaLoubes07} provide a semiparametric method for the estimation of the shifts. They use a Discrete Fourier Transform to transport the model \eqref{Sempar} into the Fourier domain. 
The important contribution of Vimond \cite{Vimond10} generalizes this study, adding the estimation of scale and height parameters. 
When the parameter $\theta$ is supposed to be random, Castillo and Loubes \cite{CastilloLoubes09} provide sharp estimators of $\theta$, following the approach of Dalalyan \textit{et al.} \cite{DGT06} in the Gaussian white noise framework. Then, they recover the unknown density of $\theta$ using a kernel density estimator. 
In this work, for the estimation of $\theta$, we propose to make use of a multidimensional version of the Robbins-Monro algorithm \cite{RobbinsMonro51}.
Assume that one can find a function $\phi$ : $\mathbb{R}^{p}\rightarrow{\mathbb{R}^{p}}$, free of the parameter $\theta$, such that $\phi(\theta)=0$. Then, it is possible to estimate
$\theta$ by the Robbins-Monro algorithm

\vspace*{-2mm}
\begin{equation}
\label{RMalgo}
\widehat{\theta}_{n+1}=\widehat{\theta}_{n}+\gamma_{n}T_{n+1}
\end{equation}

\noindent where $(\gamma_n)$ is a positive sequence of real numbers decreasing towards zero and
$(T_n)$ is a sequence of random vectors such that 
$\dE[T_{n+1}|\cF_n]=\phi(\widehat{\theta}_{n})$
where $\mathcal{F}_{n}$ stands for the $\sigma $-algebra of the events occurring up to time $n$.
Under standard conditions on the function $\phi$ and on the sequence $(\gamma_n)$, it is well-known \cite{Duflo97} that
$\widehat{\theta}_{n}$ tends to $\theta$ almost surely. The asymptotic normality of $\widehat{\theta}_{n}$ may be found in \cite{Pelletier298} whereas the quadratic strong law and the law of iterated logarithm are established in \cite{Pelletier98}. Results for randomly truncated
version of the Robbins-Monro algorithm are given in \cite{KushnerYin03}.
\par
Our second goal concerns the estimation of the scale parameter $a$. The estimation of the scale parameters, and more particularly of the sign of the scale parameters, is very important as we can see in the numerical illustrations on daily average temperatures in \cite{Vimond10}. Here, we obtain a strongly consistent estimate of the scale parameters $a$ using the prior estimation of $\theta$.
\par
The last main theoretical part of the paper is referred to the nonparametric estimation of the common shape function $f$.
A wide range of literature is available on nonparametric estimation of a regression function. We refer the reader to
\cite{Devroye01},  \cite{Tsybakov04} for two excellent books on density and regression function estimation.
Here, we focus our attention on the Nadaraya-Watson estimator \cite{Nadaraya64} \cite{Watson64} of $f$. 
More precisely, we propose to make use of a weighted recursive Nadaraya-Watson estimator \cite{Duflo97} of $f$ which
takes into account the previous estimation of $v$, $\theta$ and $a$, respectively by $\wh{v}_{n}$, $\widehat{\theta}_{n}$ and $\wh{a}_{n}$.
It is given, for all $x\in \dR$, by
\begin{equation}
\label{RNW}
\widehat{f}_{n}(x)=\sum_{j=1}^{p}\omega_{j}(x)\widehat{f}_{n,j}(x)
\end{equation}
with
$$
\widehat{f}_{n,j}(x)=\frac{1}{\wh{a}_{n,j}}\frac{\sum_{i=1}^{n} W_{i,j}(x)\left(Y_{i,j}-\wh{v}_{i-1,j}\right)}{\sum_{i=1}^{n} W_{i,j}(x)}
$$
and
$$
W_{n,j}(x)=\frac{1}{h_{n}}K\Bigl(\frac{X_{n}-\widehat{\theta}_{n-1,j}-x}{h_{n}}\Bigr),
$$
where $\wh{v}_{n,j}$, $\widehat{\theta}_{n-1,j}$ and $\wh{a}_{n,j}$ are respectively the $j$-th component of $\wh{v}_{n}$, $\widehat{\theta}_{n-1}$ and $\wh{a}_{n}$. Moreover, the kernel $K$ is a chosen probability density function and  
the bandwidth $(h_n)$ is a sequence of positive real numbers  
decreasing to zero. The main difficulty arising here is that we have to deal with the additional term
$\widehat{\theta}_{n}$ inside the kernel $K$. 
\par

\vspace{2mm}
The paper is organized as follows. Section 2 presents the model and the hypothesis which are necessary to carry out our statistical analysis. Section 3 is devoted to the parametric estimation
of the vector $v$, while Section 4 deals with our Robbins-Monro procedure for the parametric estimation of $\theta$. Section 5 concerns the parametric estimation of the vector $a$. 
In these three sections, we establish the almost sure convergence of $\widehat{v}_{n}$, $\widehat{\theta}_{n}$ and $\widehat{a}_{n}$ 
as well as  their asymptotic normality. A quadratic strong law is also provided for these three estimates. Section 6 
deals with the nonparametric estimation of $f$. Under standard regularity assumptions on the kernel $K$,
we prove the almost sure pointwise convergence of $\widehat{f}_{n}$ to $f$. In addition,
we also establish the asymptotic normality of $\widehat{f}_{n}$. 
Section 7 contains some numerical experiments on simulated and real data, illustrating
the performances of our semiparametric estimation procedure. 
The proofs of the parametric results are given is Section 8, while those concerning 
the nonparametric results are postponed to Section 9.  Finally, Section 10 is devoted to the identifiability constraints associated with the model \eqref{Sempar}.


\section{MODEL AND HYPOTHESIS}


We shall focus our attention on the model given, 
for $1\leq{j}\leq{p}$, with $p\geq2$, and for all $1\leq{i}\leq{n}$, by \eqref{Sempar}
\begin{equation*}
Y_{i,j}=a_{j}f(X_{i}-\theta_{j})+v_{j}+\varepsilon_{i,j},
\end{equation*}
where, for all $1\leq{j}\leq{p}$, $a_{j}\neq0$.
For all $1\leq{i}\leq{n}$, the noise $\left(\varepsilon_{i,j}\right)$ is a sequence of independent random variables with mean zero and variances $\dE\left[\varepsilon_{i,j}^2\right]=\sigma_{j}^2$, and independent of the random points $X_{i}$. 
In addition, as in \cite{BF10}, we make the following hypothesis.
\begin{displaymath}
\begin{array}{ll}
(\mathcal{H}_1) & \textrm{The observation times $(X_i)$ are independent and identically distributed}\\
& \textrm{with probability density function $g$, positive on its support $[-1/2;1/2]$. In} \\
& \textrm{addition, $g$ is continuous on $\dR$, twice differentiable with bounded derivatives}. \\
(\mathcal{H}_2) & \textrm{The shape function $f$ is symmetric, bounded, periodic with period 1}.
 \end{array}
\end{displaymath}


Our goal is to estimate the parameters $a$, $\theta$, $v$ as well as the common shape function $f$. However, the shape invariant model \eqref{Sempar} is not
always identifiable. Indeed, for a given vector of parameters $\left(a,\theta,v\right)$ and a given shape function $f$, one can find another vector of parameters $\left(a^{*},\theta^{*},v^{*}\right)$
and an other shape function $f^{*}$ such that for all $1\leq{j}\leq{p}$ and for all $x\in{\dR}$,
\begin{equation}
\label{identifiability}
a_{j}f(x-\theta_{j})+v_{j}=a^{*}_{j}f^{*}(x-\theta^{*}_{j})+v^{*}_{j}.
\end{equation}
H{\"a}rdle and Marron \cite{HardleMarron90} or Kneip and Engel \cite{KneipEngel95} were among the first to discuss about identifiability. 
However, in most papers dealing with shape invariant models, the identifiability is not so clear. Indeed, the identifiability conditions proposed by some authors could be not as restrictive as it should be. For example, under the conditions provided by Vimond \cite{Vimond10}, it is possible to find two different vectors of parameters $\left(a,\theta,v\right)$ and $\left(a^{*},\theta^{*},v^{*}\right)$ satisfying their conditions and \eqref{identifiability}, while under the conditions provided by Wang \textit{et al.} \cite{WangBrown2003}, it is possible to find two different shape functions $f$ and $f^{*}$ satisfying their conditions and \eqref{identifiability}.
That is the reason why we have chosen carefully the identifiability constraints. More precisely, we impose the following conditions.
\begin{displaymath}
\begin{array}{lll}
(\mathcal{H}_3) & \textrm{$\displaystyle\int_{-1/2}^{1/2}f(x)dx=0$}, & \hspace{10cm} \\
(\mathcal{H}_4) & \textrm{$a_{1}=1$, $\theta_{1}=0$ and $\underset{1\leq{j}\leq{p}} \max|\theta_{j}|<1/4$}. & \hspace{10cm}
 \end{array}
\end{displaymath}
Hypothesis $(\mathcal{H}_3)$ allows us to define uniquely the $v_{j}$ and the second constraint in $(\mathcal{H}_4)$ define uniquely the $\theta_{j}$, while $a_{1}=1$, $\theta_{1}=0$ implies that the first curve is taken as a reference. These conditions are well-adapted to our framework. Note that 
$(\mathcal{H}_3)$ and $(\mathcal{H}_4)$ could be replaced by $(\mathcal{H}_3^{\prime})$ and $(\mathcal{H}_4^{\prime})$ defined as
\begin{displaymath}
\begin{array}{lll}
(\mathcal{H}_3^{\prime}) & \textrm{$\displaystyle\int_{-1/2}^{1/2}f(x)dx=0$ and $\underset{x\in{[0;1]}}\sup|f(x)|=1$}, & \hspace{10cm} \\
(\mathcal{H}_4^{\prime}) & \textrm{$\theta_{1}=0$, $\underset{1\leq{j}\leq{p}}\min a_{j}>0$ and $\underset{1\leq{j}\leq{p}} \max|\theta_{j}|<1/2$}. & \hspace{10cm}
 \end{array}
 \end{displaymath}
An other alternative could be to substitute the hypothesis $(\mathcal{H}_4)$ by $(\mathcal{H}_4^{\prime\prime})$ or $(\mathcal{H}_4^{\prime\prime\prime})$ defined as
\begin{displaymath}
\begin{array}{lll}
(\mathcal{H}_4^{\prime\prime}) & \textrm{$a_{1}=1$, $\theta_{1}=0$, $\underset{1\leq{j}\leq{p}}\min a_{j}>0$, $\underset{1\leq{j}\leq{p}} \max|\theta_{j}|<1/2$},& \hspace{10cm}\\
(\mathcal{H}_4^{\prime\prime\prime}) & \textrm{$a_{1}>0$, $\theta_{1}=0$, $\sum_{j=1}^{p}a_{j}^{2}=p$ and $\underset{1\leq{j}\leq{p}} \max|\theta_{j}|<1/4$}. & \hspace{10cm}
 \end{array}
\end{displaymath}


Finally, note that the hypothesis of symmetry of $f$ is not necessary to ensure identifiability of the model \eqref{Sempar} but this hypothesis makes the estimation of $\theta$ easier, as we shall see in Section 4 below. All those identifiability conditions are discussed in Section 10.

\section{ESTIMATION OF THE HEIGHT PARAMETERS}

Via $(\mathcal{H}_2)$ and $(\mathcal{H}_3)$, it is not difficult to see that 
$$
\dE\left[\frac{Y_{i,j}}{g(X_{i})}\right]=v_{j}.
$$
Then, a natural consistent estimator $\widehat{v}_{n}$ of $v$ is given, for all $1\leq{j}\leq{p}$, by
\begin{equation}
\label{estimv}
\widehat{v}_{n,j}=\frac{1}{n}\sum_{i=1}^{n}\frac{Y_{i,j}}{g(X_{i})}.
\end{equation}
In order to establish the asymptotic behaviour of $\widehat{v}_{n}$, we denote by $Y$ the vector
\begin{equation}
\label{defY}
Y=\begin{pmatrix}Y_{1}\\
\vdots\\
Y_{p}\end{pmatrix},
\end{equation}
where
$$
Y_{j}=a_{j}f(X-\theta_{j})+v_{j}+\varepsilon_{j},
$$
with $X$ a random variable sharing the same distribution as the sequence $(X_{n})$ and for $1\leq{j}\leq{p}$, $\varepsilon_{j}$ sharing the same
distribution as the sequence $(\varepsilon_{i,j})$.\\
\\
The asymptotic results for $\widehat{v}_{n}$ are as follows.
\begin{thm}
\label{cvexternalshift}
Assume that $(\mathcal{H}_{1})$ to $(\mathcal{H}_{4})$ hold. Then, we have the a.s. convergence
\begin{equation}
\label{cvpsv}
\underset{n\rightarrow{+\infty}}\lim\widehat{v}_{n}=v\hspace{6mm}\textnormal{a.s.}
\end{equation}
and the asymptotic normality
\begin{equation}
\label{tlcv}
\sqrt{n}\left(\widehat{v}_{n}-v\right)\liml \mathcal{N}_{p}\left(0,\Gamma(v) \right),
\end{equation}
where $\Gamma(v)$ stands for the covariance matrix given by
 $$\Gamma(v)=\textnormal{Cov}\left(\frac{Y}{g(X)}\right).$$
 In addition, we also have the quadratic strong law
 \begin{equation}
\label{lfqv}
\underset{n\rightarrow{+\infty}}\lim\frac{1}{\log(n)}\sum_{i=1}^{n}\left(\widehat{v}_{i}-v\right)\left(\widehat{v}_{i}-v\right)^{T}=\Gamma(v)\hspace{6mm}\textnormal{a.s.}
\end{equation}
\end{thm}


\section{ESTIMATION OF THE TRANSLATION PARAMETERS}


In all the sequel, we introduce an auxiliary function $\phi$ defined for all $t\in \dR^p$, by
\begin{equation}
\label{defphi}
 \phi(t)=\dE\left[D(X,t)\begin{pmatrix}a_{1}f(X-\theta_{1})\\
\vdots\\
a_{p}f(X-\theta_{p})\end{pmatrix} \right] 
\end{equation}
where $D(X,t)$ stands for the diagonal square matrix of order $p$ defined by
\begin{equation}
\label{matrice}
D(X,t) = \frac{1}{g(X)}\mathrm{diag}\Big(\sin(2\pi(X-t_{1})),\dots,\sin(2\pi(X-t_{p}))\Big).
\end{equation}
Using the symmetry of $f$, the same calculations as in \cite{BF10} lead, for all $1\leq{j}\leq{p}$, to
\vspace{1ex}
\begin{equation}
\label{DEFPHI}
\dE\left[\frac{\sin(2\pi(X-t_{j})) }{g(X)}a_{j}f(X-\theta_{j})\right]=a_{j}f_{1} \sin(2\pi(\theta_{j}-t_{j}))
\vspace{1ex}
\end{equation} 
 where $f_1$ stands for the first Fourier coefficient of $f$
$$ f_1=\int_{-1/2}^{1/2}\cos(2\pi x)f(x)\,dx. $$
Consequently,
\begin{equation}
\label{defphi2}
 \phi(t)=f_{1} \begin{pmatrix}a_{1}\sin(2\pi(\theta_{1}-t_{1}))\\
\vdots\\
a_{p}\sin(2\pi(\theta_{p}-t_{p}))\end{pmatrix}.
\end{equation}
From now and for all the following, in order to avoid tedious discussion, we assume that $f_1 \neq 0$.
Roughly speaking, we are going to implement our Robbins-Monro procedure as in \cite{BF10} for each component of $\theta$. More precisely, for all $1\leq{j}\leq{p}$,
by noting $\phi_{j}(t)=a_{j}f_{1} \sin(2\pi(\theta_{j}-t_{j}))$, if $|t_{j}-\theta_{j}|<1/2$, $(t_{j}-\theta_{j})\phi_{j}(t)<0$ if $\textnormal{sign}(a_{j}f_{1})>0$ and $(t_{j}-\theta_{j})\phi_{j}(t)>0$ otherwise.
 Moreover, denote $K=\left[-1/4;1/4\right]$.
Then, we define the projection of $x\in{\dR}$ on $K$ by

  \begin{eqnarray*}
   \pi_{K}(x) = \left \{ \begin{array}{lll}
   \,\, \,x & \ \text{ if } \ |x|\leq 1/4, \vspace{1ex} \\
    \,1/4 & \ \text{ if } \ x\geq 1/4,  \vspace{1ex} \\
    \!\!-1/4 & \ \text{ if } \ x\leq -1/4. \\
   \end{array} \nonumber \right.
\end{eqnarray*}

\noindent Let $(\gamma_n)$ be a decreasing sequence of positive real numbers satisfying
\begin{equation}
\label{hypgamma}
\sum_{n=1}^\infty\gamma_{n}=+\infty
\hspace{1cm}\text{and}\hspace{1cm}
\sum_{n=1}^\infty\gamma_{n}^2<+\infty.
\end{equation}
For the sake of clarity, we shall make use of $\gamma_n=1/n$.
Then, for $1\leq{j}\leq{p}$, we estimate $\theta_{j}$ via the sequence $(\wh{\theta}_{n,j})$ defined, for all $n\geq{1}$, by
\begin{equation}
\label{algoRM}
\wh{\theta}_{n+1,j}=\pi_{K}\Bigl(\wh{\theta}_{n,j}+\gamma_{n+1}\text{sign}\left(a_{j}f_{1}\right)T_{n+1,j}\Bigr)
\end{equation}

\noindent where the initial value $\wh{\theta}_{0}\in K^{p}$ and the random vector $T_{n+1}$ is given by
\begin{equation}
\label{DefT}
T_{n+1}=D(X_{n+1},\wh{\theta}_{n})\begin{pmatrix}Y_{n+1,1}\\
\vdots\\
Y_{n+1,p}\end{pmatrix}.
\end{equation}

\noindent The almost sure convergence for the estimator $\wh{\theta}_{n}$ is as follows.

\begin{thm}
\label{cvpstheta}
Assume that $(\mathcal{H}_{1})$ to $(\mathcal{H}_{4})$ hold. Then, 
$\wh{\theta}_{n}$ converges almost surely to $\theta$ as $n$ goes to $+\infty$.
In addition, for $1\leq{j}\leq{p}$, the number of times that the random variable $\wh{\theta}_{n,j}+\gamma_{n+1}\text{sign}\left(a_{j}f_{1}\right)T_{n+1,j}$ goes outside of $K$ 
is almost surely finite.
\end{thm}

\begin{rem}
\label{remarquetheta}
At first sight, the estimation procedure needs the knowledge of the sign of $a_{j}f_{1}$ for all $1\leq{j}\leq{p}$. However, it is possible to do without. Indeed, denote by $(\wh{\theta}_{n,j}^{+})$ the sequence
defined, for $n\geq1$, as
$$
\wh{\theta}_{n+1,j}^{+}=\pi_{K}\Bigl(\wh{\theta}_{n,j}^{+}+\gamma_{n+1}T_{n+1,j}^{+}\Bigr),
$$
where 
\begin{equation*}
T_{n+1}^{+}=D(X_{n+1},\wh{\theta}_{n}^{+})\begin{pmatrix}Y_{n+1,1}\\
\vdots\\
Y_{n+1,p}\end{pmatrix},
\end{equation*}
and by $\wh{\theta}_{n,j}^{-}$ the sequence
defined, for $n\geq1$, as
$$
\wh{\theta}_{n+1,j}^{-}=\pi_{K}\Bigl(\wh{\theta}_{n,j}^{-}-\gamma_{n+1}T_{n+1,j}^{-}\Bigr),
$$
where
\begin{equation*}
T_{n+1}^{-}=D(X_{n+1},\wh{\theta}_{n}^{-})\begin{pmatrix}Y_{n+1,1}\\
\vdots\\
Y_{n+1,p}\end{pmatrix}.
\end{equation*}
Then, two events are possible. More precisely, for $1\leq{j}\leq{p}$,
$$
\underset{n\rightarrow{+\infty}}\lim \wh{\theta}_{n,j}^{+}=\theta_{j}\hspace{8mm}\textnormal{and}\hspace{8mm}\underset{n\rightarrow{+\infty}}\lim |\wh{\theta}_{n,j}^{-}|=1/4\hspace{5mm}\textnormal{a.s.}
$$
or
$$
\underset{n\rightarrow{+\infty}}\lim \wh{\theta}_{n,j}^{-}=\theta_{j}\hspace{8mm}\textnormal{and}\hspace{8mm}\underset{n\rightarrow{+\infty}}\lim |\wh{\theta}_{n,j}^{+}|=1/4\hspace{5mm}\textnormal{a.s.}
$$
Hence, for $n$ large enough, the vector $\wh{\theta}_{n}$ which is considered is the vector whose absolute value of the j-th component is given by 
$\min(|\wh{\theta}_{n,j}^{+}|,|\wh{\theta}_{n,j}^{-}|).$
Nevertheless, for the sake of clarity, we shall do as if the sign of $a_{j}f_{1}$ is known.
\end{rem}
\vspace{1ex}

In order to establish the asymptotic normality of $\wh{\theta}_{n}$, it is necessary to
introduce a second auxiliary function $\varphi$ defined, for all $t\in \dR^p$, as

\begin{equation}
 \label{defvarphi}
\varphi(t)=\mathbb{E}[V(t)V(t)^{T}]
\end{equation}
where $V(t)$ is given by
$$
V(t)=\mathrm{diag}\Big(\textnormal{sign}(a_{1}f_{1}),\dots,\textnormal{sign}(a_{p}f_{1})\Big)D(X,t)Y,
$$
with $Y$ given by \eqref{defY}. As soon as $4\pi |f_{1}|\underset{1\leq j\leq p} \min|a_{j}|>1$, denote for all $1\leq{k,l}\leq{p}$,
$$
\Sigma(\theta)_{k,l}=\frac{\varphi(\theta)_{k,l}}{2\pi( |a_{k}|+|a_{l}|)|f_{1}|-1}.
$$
\begin{thm}
\label{thmcltrm}
Assume that $(\mathcal{H}_{1})$ to $(\mathcal{H}_{4})$ hold. In addition, suppose that $(\veps_{i,j})$ has a finite moment of order $>2$ and that
$$4\pi |f_{1}|\underset{1\leq j\leq p} \min|a_{j}|>1.$$
Then, we have the asymptotic normality
\begin{equation}
\label{cltrm}
\sqrt{n}(\wh{\theta}_{n}-\theta) \liml \cN_{p}(0, \Sigma(\theta)).
\end{equation}
\end{thm}

\begin{thm}
\label{thmlilqsl}
Assume that $(\mathcal{H}_{1})$ to $(\mathcal{H}_{4})$ hold. In addition, suppose that $(\veps_{i,j})$ has a finite moment of order $>2$ and that
$$4\pi |f_{1}|\underset{1\leq j\leq p} \min|a_{j}|>1.$$ 
Then,  we have the law of iterated logarithm, given, for all $w\in{\dR^p}$, by
\begin{eqnarray}  
\label{lilrm}
\limsup_{n \rightarrow \infty} \left(\frac{n}{2 \log \log n} \right)^{1/2}
w^{T}(\wh{\theta}_{n}-\theta)
&=& - \liminf_{n \rightarrow \infty}
\left(\frac{n}{2 \log\log n}\right)^{1/2} 
w^{T}(\wh{\theta}_{n}-\theta)  \notag \\
&=& \sqrt{w^{T}\Sigma(\theta)w}
\hspace{1cm}\textnormal{a.s.}
\end{eqnarray}
In particular, 
\begin{equation}
\label{lilsuprm}  
\limsup_{n \rightarrow \infty} \left(\frac{n}{2 \log \log n} \right)
\left(\wh{\theta}_{n}-\theta\right)\left(\wh{\theta}_{n}-\theta\right)^{T}=\Sigma(\theta)
\hspace{1cm}\textnormal{a.s.}
\end{equation}
In addition, we also have the quadratic strong law 
\begin{equation}
\label{qslrm}
\lim_{n\rightarrow \infty}\frac{1}{\log n}\sum_{i=1}^{n}(\wh{\theta}_{i}-\theta)(\wh{\theta}_{i}-\theta)^{T}=\Sigma(\theta)
\hspace{1cm}\textnormal{a.s.}
\end{equation}
\end{thm}

\begin{rem}
In the particular case where $$4\pi |f_{1}|\underset{1\leq j\leq p} \min|a_{j}|=1,$$ it is also possible to show, thanks to Theorem 2.2.12 page 52 of \cite{Duflo97}, that
$$
\sqrt{\frac{n}{\log(n)}}(\wh{\theta}_{n}-\theta) \liml \cN(0, \varphi(\theta)).
$$
Asymptotic results can also be established when 
$$0<4\pi |f_{1}|\underset{1\leq j\leq p} \min|a_{j}|<1.$$
However, we have chosen to focus our attention on the more attractive case
$$4\pi |f_{1}|\underset{1\leq j\leq p} \min|a_{j}|>1.$$
\end{rem}

\begin{rem}
\label{remefficiency}
We clearly have via differentiation
$$
D\phi(t)=-2 \pi f_1\mathrm{diag}\Big(a_{1}\cos(2\pi(\theta_{1} -t_{1})),\dots,
a_{p}\cos(2\pi(\theta_{p} -t_{p}))\Big).
$$
Consequently, the value $$D\phi(\theta)=-2\pi f_1\mathrm{diag}\Big(a_{1},\dots,
a_{p}\Big)$$ does not depend upon the unknown parameter $\theta$. 
On the one hand, if the first Fourier coefficient $f_1$ of $f$ and the vector $a$ are known,
it is possible to provide, via a slight modification of \eqref{algoRM}, an asymptotically efficient estimator $\wh{\theta}_{n}$ of $\theta$.
More precisely, for all $1\leq{j}\leq{p}$, it is only necessary to replace $\gamma_n=1/n$ in \eqref{algoRM} by $\gamma_n=\gamma_{j}/n$ where
$$
\gamma_{j}= \frac{1}{2 \pi a_{j}f_1}.
$$
Then, we deduce from the part 10.2.2 of the book of Kushner and Yin \cite{KushnerYin03} page 331 that $\wh{\theta}_{n}$ is an asymptotically efficient estimator of $\theta$ 
with
\begin{equation}
\label{cltrmefficient}
\sqrt{n}(\wh{\theta}_{n}-\theta) \liml \cN\Bigl(0,\ell(\theta)\Bigr).
\end{equation}
where for all $1\leq{k,l}\leq{p}$,
$$
\ell(\theta)_{k,l}=\frac{\varphi(\theta)_{k,l}}{4\pi^{2}|f_{1}|^{2}|a_{k}a_{l}|},
$$
with $\varphi(\theta)$ given by \eqref{defvarphi}.
On the other hand, if $f_1$ and $a$ are unknown, it is also possible to provide by the same procedure
an asymptotically efficient estimator $\wh{\theta}_{n}$ of $\theta$ replacing $f_1$ by its natural estimate
\begin{equation}
\label{estimf1}
\widehat{f}_{1,n}=\frac{1}{n}\sum_{i=1}^{n}\frac{\cos(2\pi X_{i})}{g(X_{i})}Y_{i,1},
\end{equation}
and replacing the value $a_{j}$ by the estimate $\widetilde{a}_{n,j}$ which is defined in the next Section.
\end{rem}

\begin{rem}
\label{nonsym}
As in \cite{BF10}, it is also possible to get rid of the symmetry assumption on $f$.
However, it requires the knowledge of the first Fourier coefficients of $f$
$$
f_1=\int_{-1/2}^{1/2}\cos(2\pi x)f(x)\,dx
\hspace{1cm}\textnormal{and}\hspace{1cm}
g_1=\int_{-1/2}^{1/2}\sin(2\pi x)f(x)\,dx.
$$
On the one hand, it is necessary to assume that $f_1\neq 0$ or $g_1\neq 0$, and to replace in \eqref{defphi} the diagonal
matrix $D(X,t)$ defined in \eqref{matrice}, by
\begin{equation}
\label{matrice2}
\Delta(X,t)=\frac{1}{g(X)}\mathrm{diag}\Big(\delta(X,t_{1}),\dots,\delta(X,t_{p})\Big),
\end{equation}
where, for all $1\leq{j}\leq{p}$,
$$
\delta(X,t_{j})=f_{1}\sin(2\pi(X-t_{j}))-g_{1}\cos(2\pi(X-t_{j})).
$$
Then, Theorem \ref{cvpstheta} is true for the projected Robbins-Monro algorithm defined, for all $1\leq{j}\leq{p}$, by
\begin{equation*}
\wh{\theta}_{n+1,j}=\pi_{K}\Bigl(\wh{\theta}_{n,j}+\textnormal{sign}(a_{j})\gamma_{n+1}T_{n+1,j}\Bigr),
\end{equation*}
where the initial value $\wh{\theta}_{0} \in K^{p}$ and the random vector $T_{n+1}$ is given by
\begin{equation*}
T_{n+1}=\Delta(X_{n+1},\wh{\theta}_{n})\begin{pmatrix}Y_{n+1,1}\\
\vdots\\
Y_{n+1,p}\end{pmatrix}.
\end{equation*}
On the other hand, we also have to replace the second function $\varphi$ defined in \eqref{defvarphi} by
$$
\Psi(t)=\mathbb{E}[W(t)W(t)^{T}],
$$
where $W(t)$ is given by
$$
W(t)=\mathrm{diag}\Big(\textnormal{sign}(a_{1}f_{1}),\dots,\textnormal{sign}(a_{p}f_{1})\Big)\Delta(X,t)Y$$
Then, as soon as $4\pi (f_{1}^2+g_{1}^2)\underset{1\leq{j}\leq{p}}\min|a_{j}|>1$, Theorems \ref{thmcltrm} and \ref{thmlilqsl} hold with $\Sigma(\theta)$ given for all $1\leq{k,l}\leq{p},$ by
$$
\Sigma(\theta)_{k,l}=\frac{\Psi(\theta)_{k,l}}{2\pi (f_{1}^2+g_{1}^2)(|a_{k}|+|a_{l}|) -1}.
$$
In the rest of the paper, we shall not go in that direction as our strategy is to make very few assumptions on the Fourier coefficients
of $f$.
\end{rem}


\section{ESTIMATION OF THE SCALE PARAMETERS}


Henceforth, we introduce an other auxiliary function $\psi$ defined for all $t\in \dR^p$, by
\begin{equation}
\label{defpsi}
 \psi(t)=\dE\left[C(X,t)\begin{pmatrix}a_{1}f(X-\theta_{1})\\
\vdots\\
a_{p}f(X-\theta_{p})\end{pmatrix} \right] 
\end{equation}
where $C(X,t)$ is the diagonal matrix of order $p$, given by 
\begin{equation}
\label{matrice3}
C(X,t)=\frac{1}{g(X)}\mathrm{diag}\Big(\cos(2\pi(X-t_{1})),\dots,\cos(2\pi(X-t_{p}))\Big).
\end{equation}
\vspace{2ex}
\noindent As for \eqref{defphi2}, we have
\vspace{1ex}
\begin{equation}
\label{defpsi2}
 \psi(t)=f_{1}\begin{pmatrix}a_{1}\cos(2\pi(\theta_{1}-t_{1}))\\
\vdots\\
a_{p}\cos(2\pi(\theta_{p}-t_{p}))\end{pmatrix}.
\end{equation}

\noindent Then, it is clear from Theorem \ref{cvpstheta} that $\psi(\widehat{\theta}_{n})$ tends to $\psi(\theta)=f_{1}a$. Hence, denote by $(\widehat{a}_{n})$ and $(\widetilde{a}_{n})$
the sequences defined, for
$n\geq1$ and for all $1\leq{j}\leq{p}$, by
\begin{equation}
\label{defa}
\widehat{a}_{n,j}=\frac{1}{n f_{1}}\sum_{i=1}^{n}\frac{\cos(2\pi(X_{i}-\widehat{\theta}_{i-1,j}))}{g(X_{i})}Y_{i,j}
\end{equation}
and
\begin{equation}
\label{defatilde}
\widetilde{a}_{n,j}=\frac{1}{n\wh{f}_{1,n}}\sum_{i=1}^{n}\frac{\cos(2\pi(X_{i}-\widehat{\theta}_{i-1,j}))}{g(X_{i})}Y_{i,j}.
\end{equation}
where $\wh{f}_{1,n}$ is given by \eqref{estimf1}. Let $I_{p}$ the identity matrix of order $p$, $e_{1}$ the first Euclidian vector of $\dR^{p}$ and $M_{p}$ the square matrix given by
\begin{equation}
\label{defMp}
M_{p}=I_{p}-ae_{1}^{T}.
\end{equation}
Then, the asymptotic behaviours of $\wh{a}_{n}$ and of $\widetilde{a}_{n}$ are as follows.

\begin{thm}
\label{cva}
Assume that $(\mathcal{H}_{1})$ to $(\mathcal{H}_{4})$ hold. Then, we have
\begin{equation}
\label{cvpsa}
\underset{n\rightarrow{+\infty}}\lim\widehat{a}_{n}=a\hspace{6mm}\textnormal{a.s.}
\end{equation}
and
\begin{equation}
\label{cvpsatilde}
\underset{n\rightarrow{+\infty}}\lim\widetilde{a}_{n}=a\hspace{6mm}\textnormal{a.s.}
\end{equation}
and the asymptotic normalities
\begin{equation}
\label{tlca}
\sqrt{n}\left(\widehat{a}_{n}-a\right)\liml \mathcal{N}_{p}\left(0,\Gamma(a) \right),
\end{equation}
and
\begin{equation}
\label{tlcatilde}
\sqrt{n}\left(\widetilde{a}_{n}-a\right)\liml \mathcal{N}_{p}\left(0,M_{p}\Gamma(a)M_{p}^{T}\right),
\end{equation}
where $\Gamma(a)$ stands for the covariance matrix given by
\begin{equation}
\label{defgammaa}
\Gamma(a)=\frac{1}{f_{1}^{2}}\textnormal{Cov}\left(C(X,\theta)Y\right).
\end{equation}
In addition, we also have the quadratic strong laws
\begin{equation}
\label{lfqa}
\underset{n\rightarrow{+\infty}}\lim\frac{1}{\log(n)}\sum_{i=1}^{n}\left(\widehat{a}_{i}-a\right)\left(\widehat{a}_{i}-a\right)^{T}=\Gamma(a)\hspace{6mm}\textnormal{a.s.}
\end{equation}
\begin{equation}
\label{lfqatilde}
\underset{n\rightarrow{+\infty}}\lim\frac{1}{\log(n)}\sum_{i=1}^{n}\left(\widetilde{a}_{i}-a\right)\left(\widetilde{a}_{i}-a\right)^{T}=M_{p}\Gamma(a)M_{p}^{T}\hspace{6mm}\textnormal{a.s.}
\end{equation}
\end{thm}

\section{ESTIMATION OF THE REGRESSION FUNCTION}

\vspace{2ex}
In this section, we are interested in the nonparametric estimation of the regression function $f$ via a recursive Nadaraya-Watson
estimator. On the one hand, we add the following standard hypothesis.
\vspace{1ex}
\begin{displaymath}
\begin{array}{ll}
(\mathcal{H}_5) & \textrm{The regression function $f$ is Lipschitz}.\hspace{9cm}
\vspace{1ex}
\end{array}
\end{displaymath}
and we suppose that $f_{1}$ is known (see Remark \ref{f1unknown} below).\\
On the other hand, we follow the same approach as in \cite{BF10}. Moreover, for more accuracy, we consider a weighted version of the Nadaraya-Watson estimator
\begin{equation}
\label{RNWS}
\wh{f}_{n}(x)=\sum_{j=1}^{p}\omega_{j}(x)\wh{f}_{n,j}(x),
\end{equation}

\vspace{-4mm}
\noindent where, for all $1\leq{j}\leq{p}$, 

\vspace{-7mm}
\begin{equation}
\label{poids}
\omega_{j}(x)=\omega_{j}(-x),\hspace{12mm}\omega_{j}(x)\geq0\hspace{12mm}\text{and}\hspace{12mm}\sum_{j=1}^{p}\omega_{j}(x)=1,
\end{equation}

\vspace{-5mm}
\begin{equation}
\label{defNWj}
\wh{f}_{n,j}(x)=\frac{1}{\widehat{a}_{n,j}}\frac{\sum_{i=1}^{n} (W_{i,j}(x)+W_{i,j}(-x))\left(Y_{i,j}-\widehat{v}_{i-1,j}\right)}{\sum_{i=1}^{n} (W_{i,j}(x)+W_{i,j}(-x))},
\end{equation}

\noindent and

\vspace{-6mm}
\begin{equation}
\label{defW}
W_{n,j}(x)=\frac{1}{h_{n}}K\Bigl(\frac{X_{n}-\widehat{\theta}_{n-1,j}-x}{h_{n}}\Bigr).
\end{equation}

\vspace{1ex}
\noindent The bandwidth $(h_n)$ is a sequence of positive real numbers, 
decreasing to zero, such that $n h_n$ tends to infinity. For the sake of simplicity, we propose to make use of
$h_n = 1/n^{\alpha}$ with $\alpha \in\, ]0,1[$.
Moreover,  we shall assume in all the sequel that the kernel $K$ is a nonnegative symmetric function, bounded with compact support, satisfying  
$$ \int_{\dR} K(x)\,dx = 1
\hspace{1cm}\textnormal{and}\hspace{1cm}
\int_{\dR} K^2(x)\,dx= \nu^2.
$$

\begin{thm}
\label{thmaspnw}
Assume that $(\mathcal{H}_{1})$ to $(\mathcal{H}_{5})$ hold and that the sequence $(\veps_{i,j})$ has a finite moment of order $>2$.
Then, for any $x\in{[-1/2;1/2]}$, we have
\begin{equation} 
\label{Cvgfchapn}
\lim_{n\rightarrow \infty}
\wh{f}_{n}(x)=f(x)\hspace{1cm}\textnormal{a.s.}
\end{equation}
\end{thm}

\begin{thm}
\label{thmcltnw}
Assume that $(\mathcal{H}_{1})$ to $(\mathcal{H}_{5})$ hold and that the sequence $(\veps_{i,j})$ has a finite moment of order $>2$. Then, as soon as
the bandwidth $(h_n)$ satisfies $h_n = 1/n^{\alpha}$ with $\alpha >1/3$, we have
for any $x\in{[-1/2;1/2]}$ with $x\neq{0}$, the pointwise asymptotic normality
\begin{equation}
\label{cltnwx}
\sqrt{nh_n}(\wh{f}_{n}(x)-f(x)) \liml \cN\left(0,\frac{\nu^2}{1+\alpha}\sum_{j=1}^{p}\frac{\sigma_{j}^{2}\omega_{j}^2(x)}{a_{j}^{2}\left(g(\theta_{j}+x)+g(\theta_{j}-x)\right)}\right).
\end{equation}
In addition, for $x=0$,
\begin{equation}
\label{cltnwzero}
\sqrt{nh_n}(\wh{f}_{n}(0)-f(0)) \liml \cN\left(0,\frac{\nu^2}{1+\alpha}\sum_{j=1}^{p}\frac{\sigma_{j}^{2}\omega_{j}^2(0)}{a_{j}^{2}g(\theta_{j})}\right).
\end{equation}
\end{thm}

\begin{rem}
\label{f1unknown}
If $f_{1}$ is unknown, it is necessary to replace the estimator $\wh{f}_{n}(x)$ defined in \eqref{RNWS} by $\widetilde{f}_{n}(x)$ given by 
\begin{equation*}
\widetilde{f}_{n}(x)=\sum_{j=1}^{p}\omega_{j}(x)\widetilde{f}_{n,j}(x),
\end{equation*}
where $\widetilde{f}_{n,j}(x)$ is defined, for all $x\in{[-1/2;1/2]}$ and for all $n\geq1$, by
\begin{equation*}
\widetilde{f}_{n,j}(x)=\frac{1}{\widetilde{a}_{n,j}}\frac{\sum_{i=1}^{n} (W_{i,j}(x)+W_{i,j}(-x))\left(Y_{i,j}-\widehat{v}_{i-1,j}\right)}{\sum_{i=1}^{n} (W_{i,j}(x)+W_{i,j}(-x))}.
\end{equation*}
One can observe that Theorems \ref{thmaspnw} and \ref{thmcltnw} are still true for $\widetilde{f}_{n}(x)$ with the same asymptotic variances.
\end{rem}

\begin{rem}
\label{rempoids}
The choice of the weights $\omega_{j}(x)$ could be important. Intuitively, the asymptotic variances in \eqref{cltnwx} and \eqref{cltnwzero} are minimal if for all $1\leq{j}\leq{p}$,
 $\omega_{j}(x)$ is inversely proportional to the variance $\sigma_{j}^{2}$ of the noise. More precisely, the Lagrange Multiplier Theorem gives us the values of the weights for
 which the asymptotic variances in \eqref{cltnwx} and \eqref{cltnwzero} are minimal under the constraint \eqref{poids}. They are given, for all $1\leq{j}\leq{p}$ and for all $x\in{[-1/2;1/2]}$, by
 $$
 \omega_{j}(x)=\frac{m_{j}(x)}{\sum_{k=1}^{p}m_{k}(x)}
 $$
 where
 $$
 m_{j}(x)=\frac{a_{j}^{2}\left(g(\theta_{j}+x)+g(\theta_{j}-x)\right)}{\sigma_{j}^{2}}.
 $$
 For these values, the asymptotic variances in \eqref{cltnwx} and \eqref{cltnwzero} are respectively given, for $x\neq{0}$, by
 $$
\frac{\nu^2}{1+\alpha}\left(\sum_{j=1}^{p}m_{j}(x)\right)^{-1},
 $$
 and for $x=0$, by
 $$
 \frac{\nu^2}{1+\alpha}\left(\sum_{j=1}^{p}\frac{m_{j}(0)}{2}\right)^{-1}.
 $$
\end{rem}


\section{SIMULATIONS}

In this section, we illustrate the asymptotic behavior of our estimates on simulated data as well as on real data.

\subsection{Simulated data}
\begin{figure}[htp]
\vspace{-2ex}
\begin{center}
\includegraphics[scale=0.53]{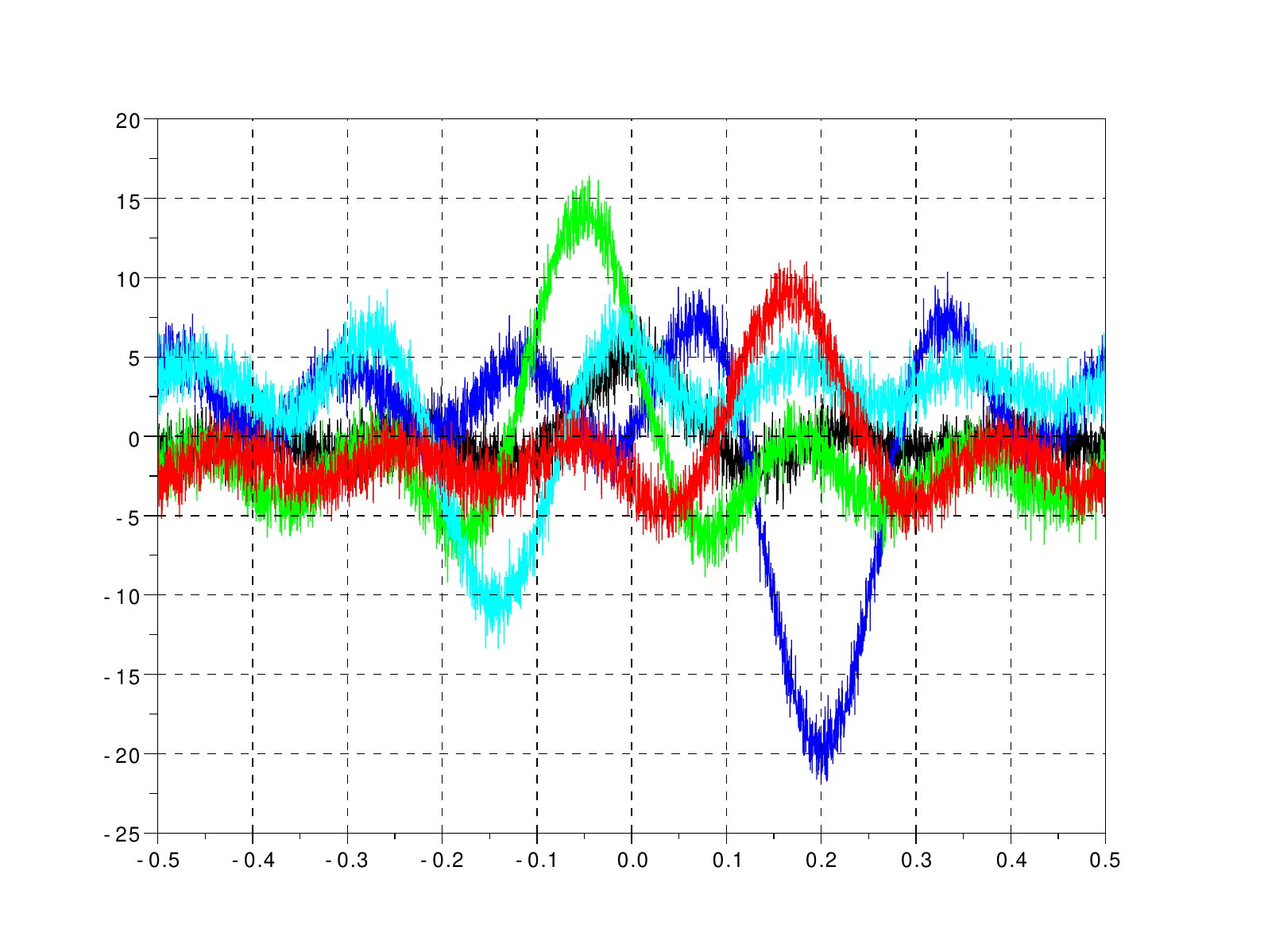}
\vspace{-2ex}
\caption{Simulated data}
\label{simdata}
\end{center}
\end{figure}

\begin{figure}[htb]
\vspace{-2ex}
\begin{minipage}[b] {0.48\linewidth}
\centering 
\centerline {\includegraphics[scale=0.4]{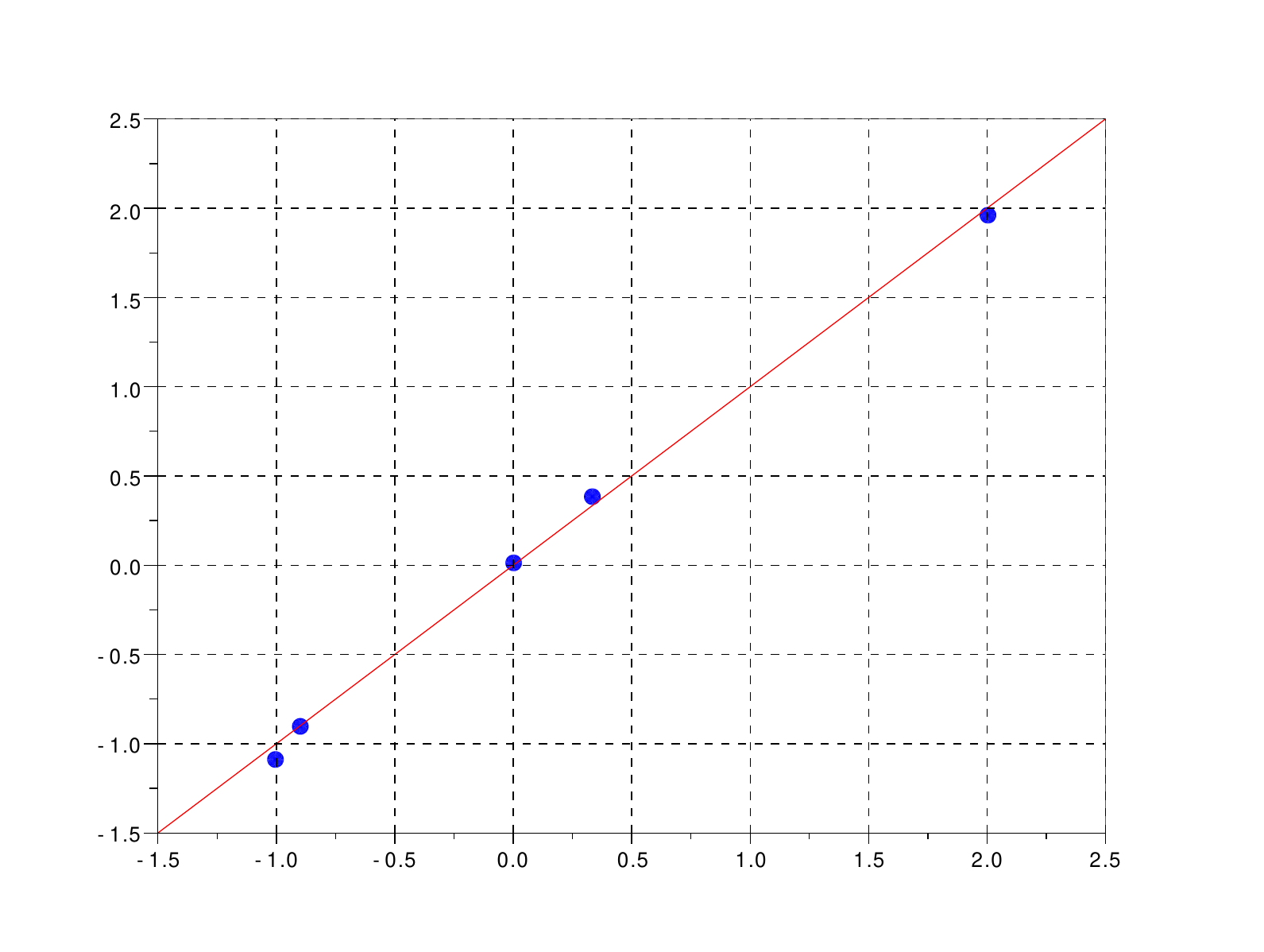} }
\vspace{0.1cm}
\medskip
\end{minipage}
\hfill
\begin{minipage}[b]{0.48\linewidth}
\centering 
\centerline {\includegraphics[scale=0.4]{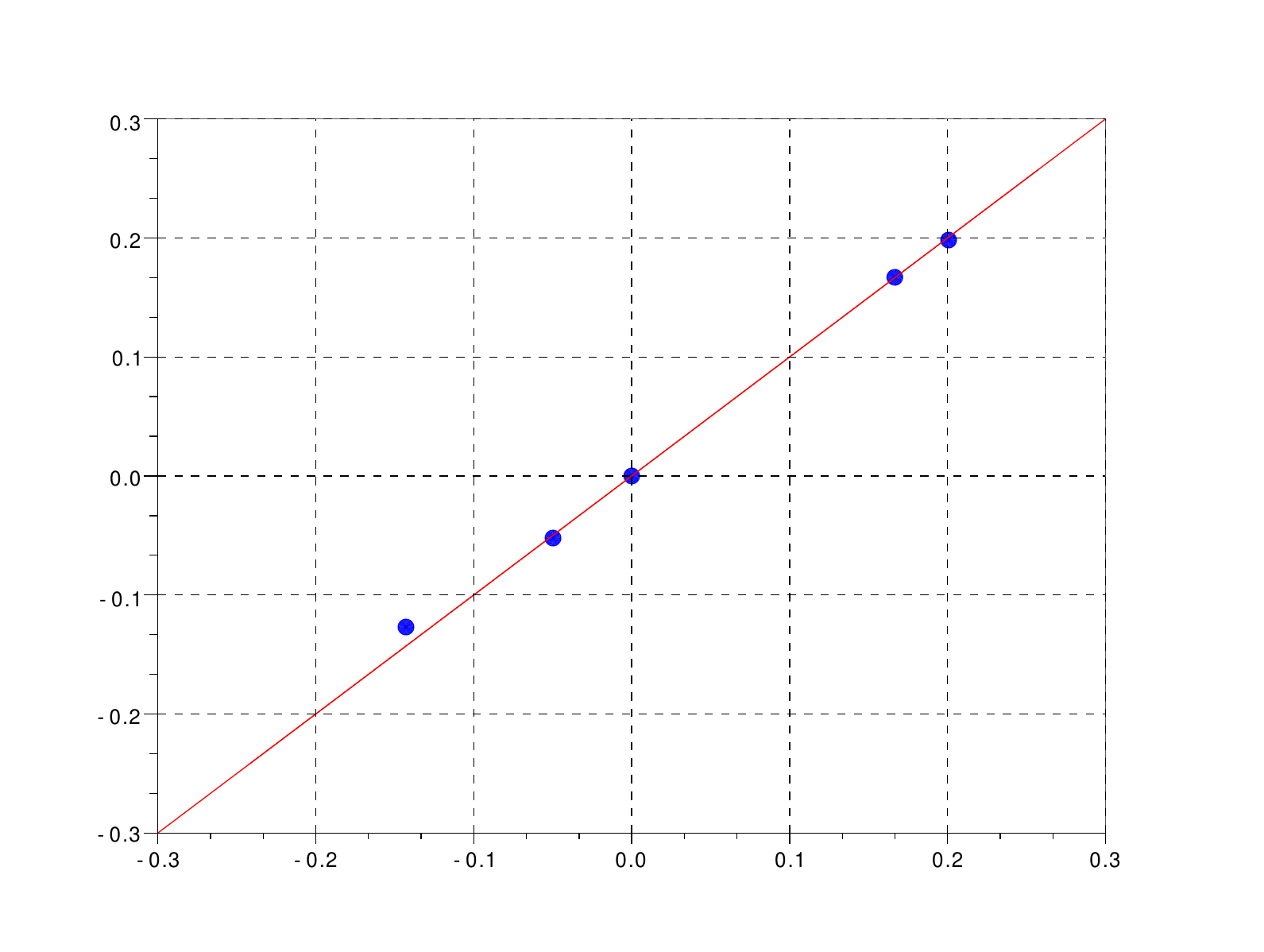} }
\vspace{0.1cm}
\medskip
\end{minipage}
\hfill
\begin{minipage}[b]{0.48\linewidth}
\centering 
\centerline {\includegraphics[scale=0.4]{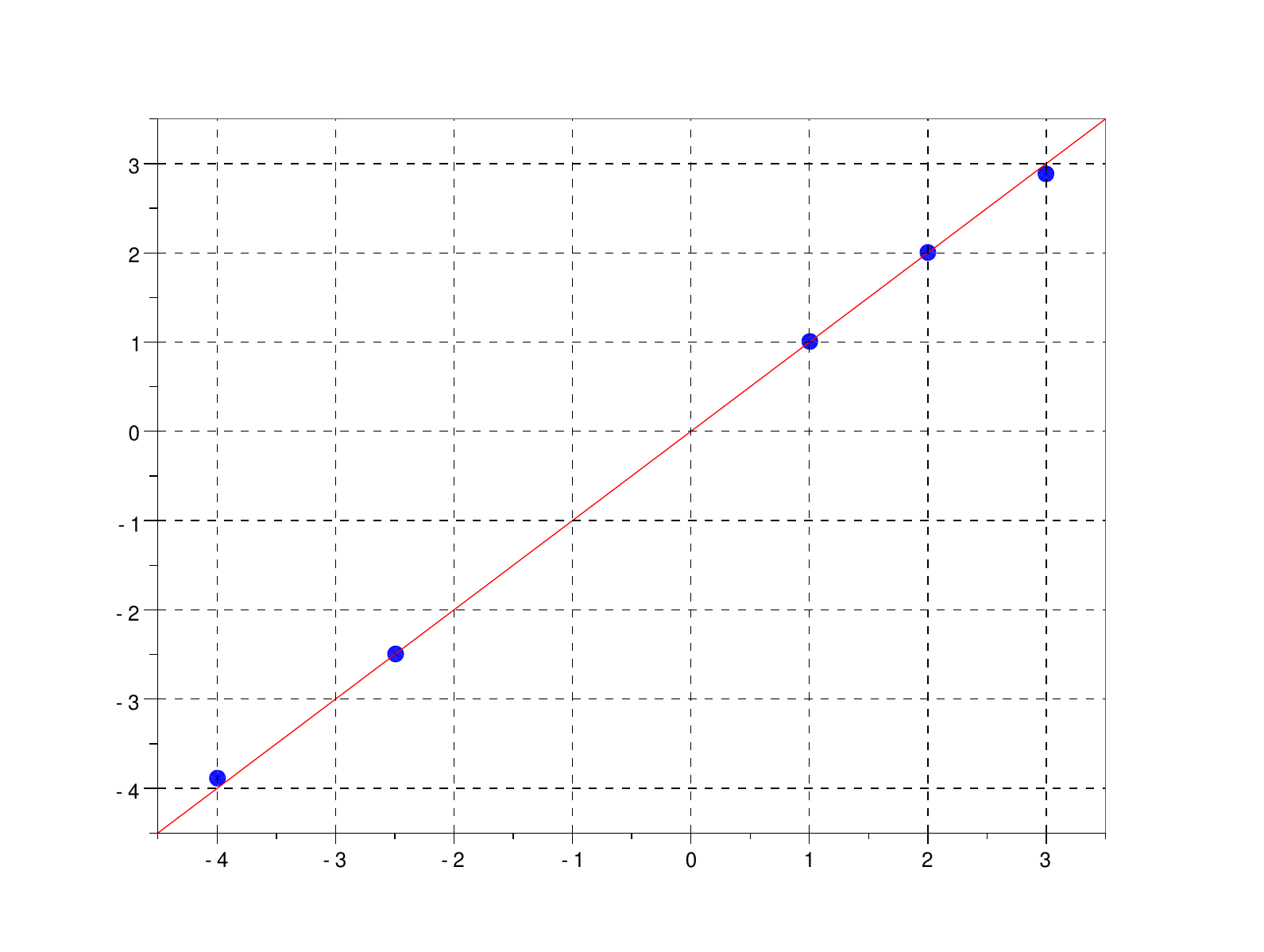} }
\vspace{0.1cm}
\medskip
\end{minipage}
\vspace{-2ex}
\caption{Estimation of $v$, $\theta$ and $a$}
\label{param}
\end{figure}

We consider data simulated according to the model \eqref{Sempar}
\begin{equation*}
Y_{i,j}=a_{j}f(X_{i}-\theta_{j})+v_{j}+\veps_{i,j}
\hspace{7mm}
\end{equation*}
where $1\leq{j}\leq{p}$ and $1\leq{i}\leq{n}$ with $p=5$ and $n=2\,000$. We have chosen the height parameters $v=\left(0,1/3,-1,2,-9/10\right)^{T}$, the translation parameters $\theta=\left(0,1/5,-1/20,-1/7,1/6\right)^{T}$ and the scale parameters $a=\left(1,-4,3,-5/2,-2\right)^{T}$. In addition, the noise $(\veps_{i,j})$ is a sequence of i.i.d. random variables with $\mathcal{N}\left(0,1\right)$ distribution. The random variables $(X_{i})$ are simulated according to the uniform distribution on $[-1/2;1/2]$ and the regression function, whose $f_{1}=1/2$, is given for all $x\in{[-1/2;1/2]}$ by 
$$
f(x)=\sum_{k=1}^{5}\cos(2k\pi x).
$$
The simulated data are given in Figure \ref{simdata}. The results of the estimates of the vectors $v$, $\theta$ and $a$ by $\wh{v}_{n}$, $\wh{\theta}_{n}$ and $\wh{a}_{n}$ are given in Figure \ref{param}. The true values are drawn in abscissa and their estimates in ordinate. One can observe that the true values and their estimates are very close, showing that our parametric estimation procedure performs pretty well on the simulated data.

\begin{figure}[htb]
\vspace{-2ex}
\begin{minipage}[b] {0.48\linewidth}
\centering 
\centerline {\includegraphics[scale=0.44]{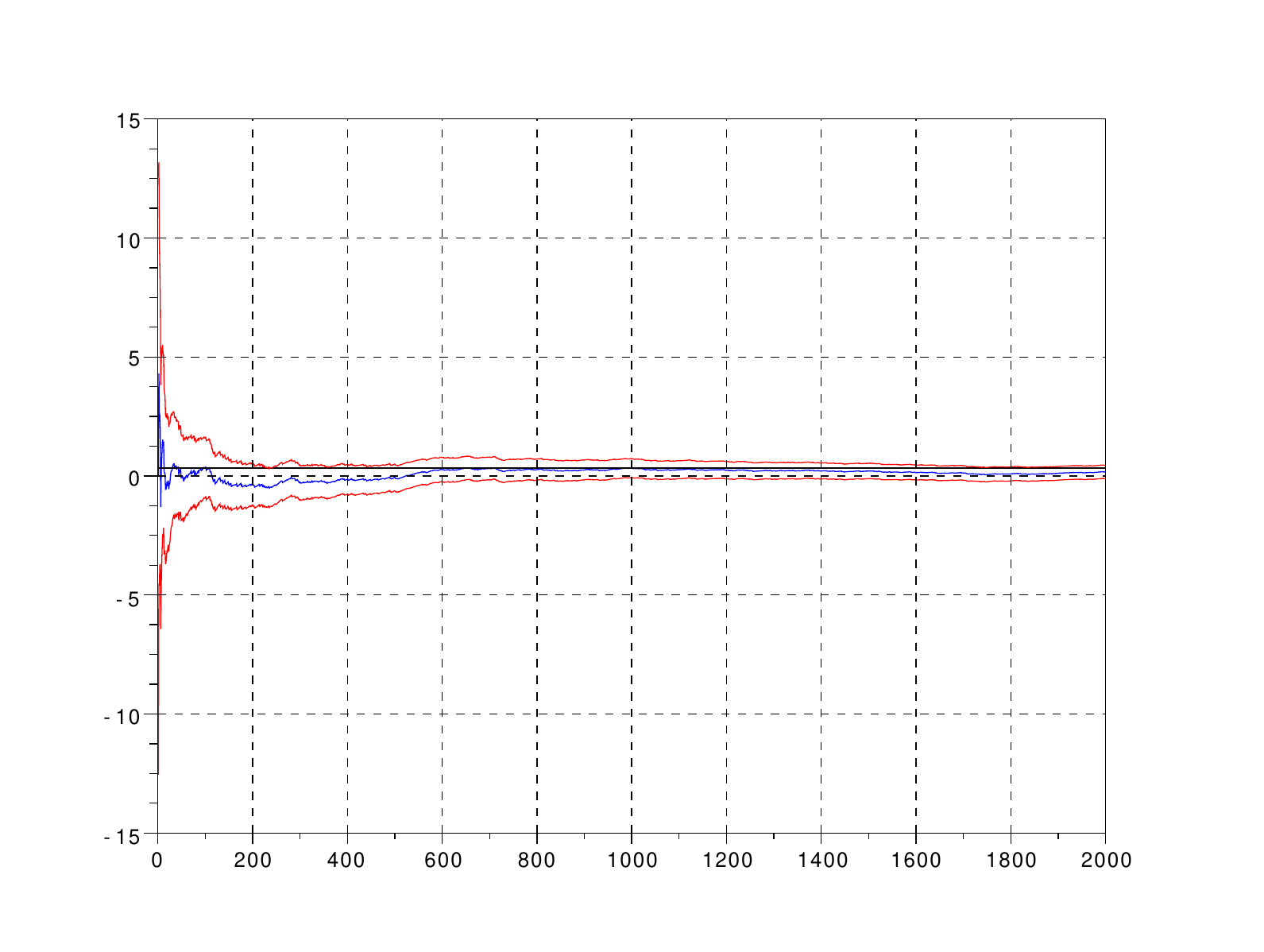} }
\vspace{0.1cm}
\medskip
\end{minipage}
\hfill
\begin{minipage}[b]{0.48\linewidth}
\centering 
\centerline {\includegraphics[scale=0.44]{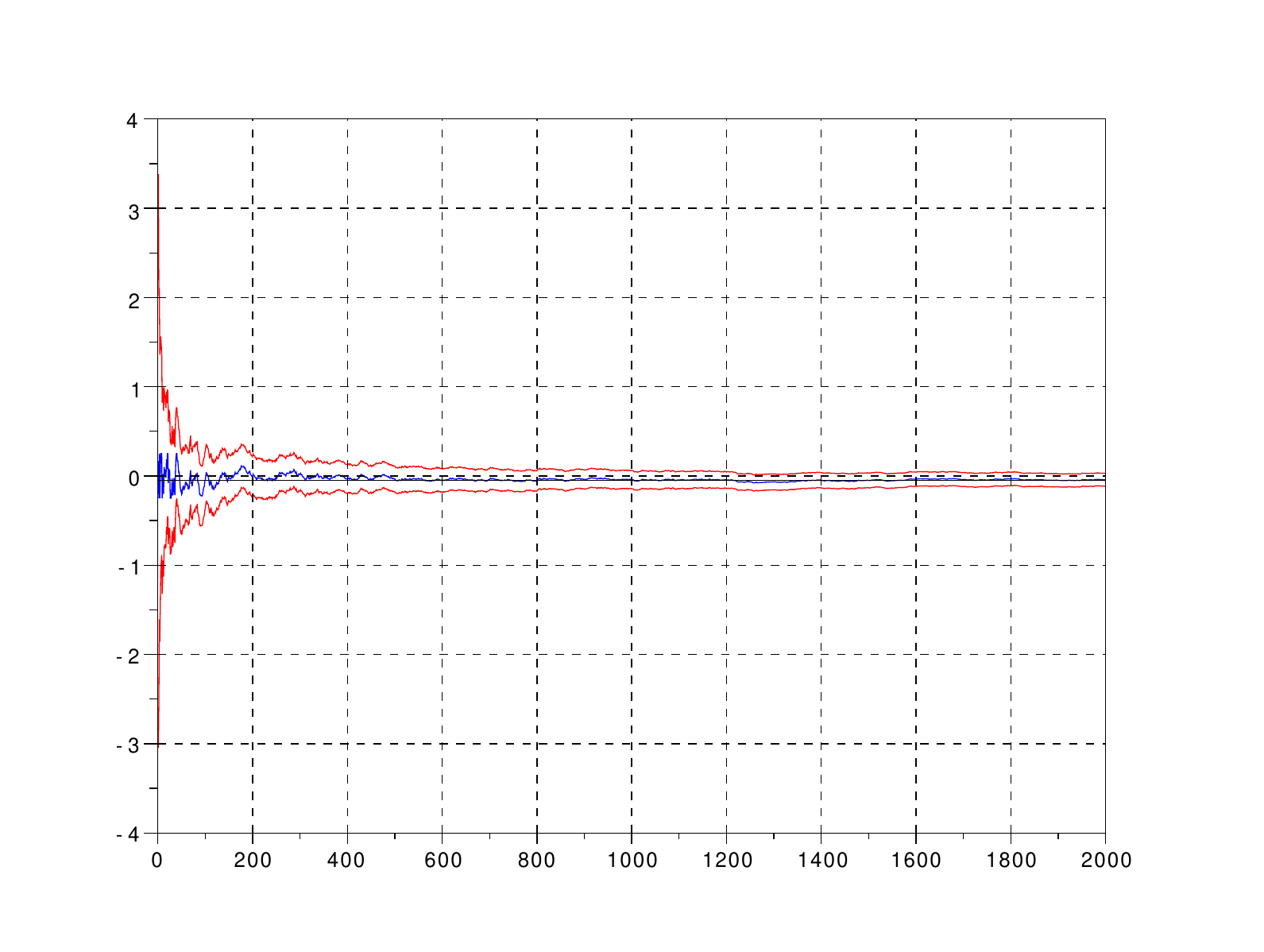} }
\vspace{0.1cm}
\medskip
\end{minipage}
\hfill
\begin{minipage}[b]{0.48\linewidth}
\centering 
\centerline {\includegraphics[scale=0.44]{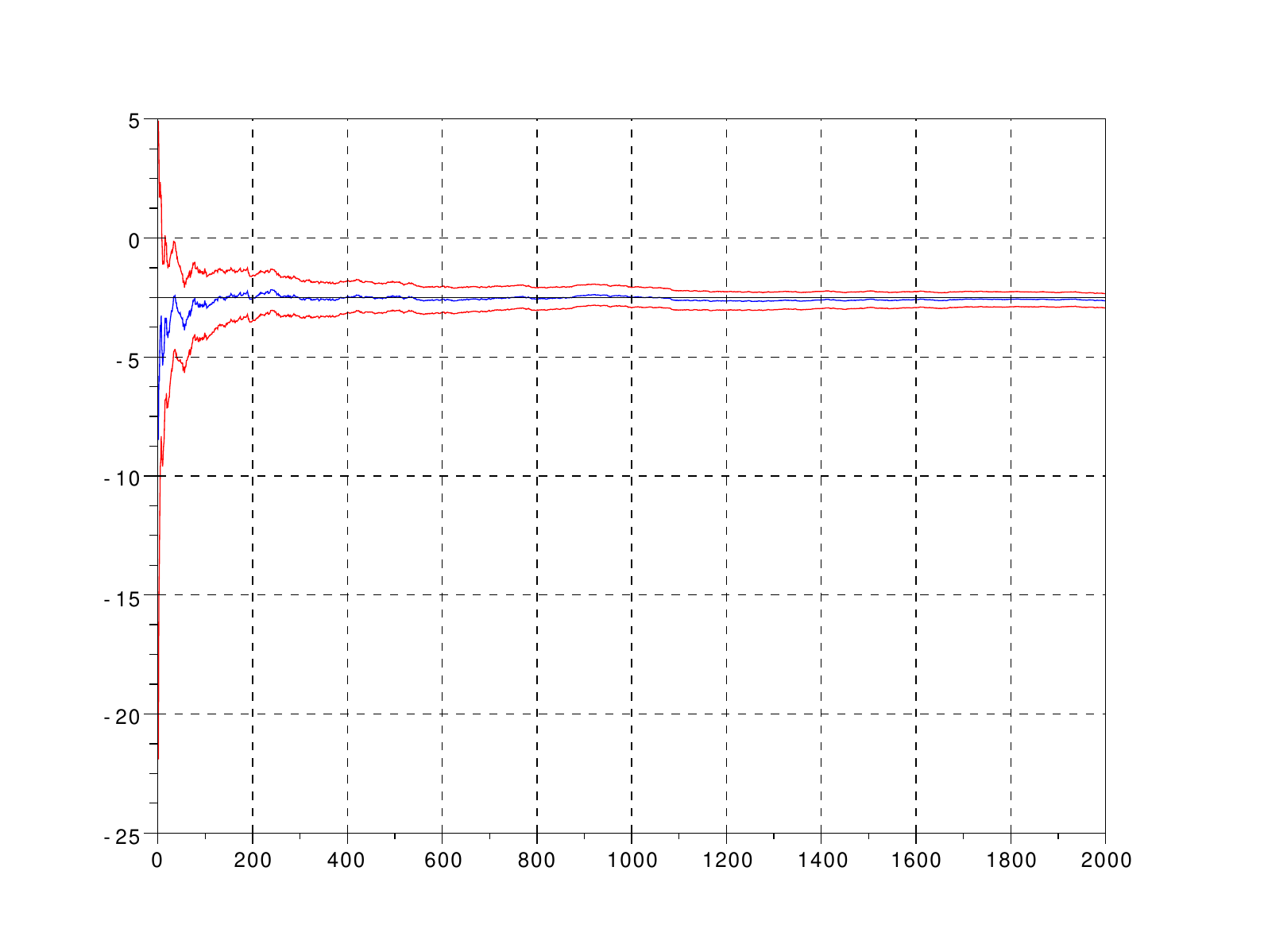} }
\vspace{0.1cm}
\medskip
\end{minipage}
\vspace{-2ex}
\caption{Confidence intervals of $v_{2}$, $\theta_{3}$ and $a_{4}$}
\label{ICparam}
\end{figure}

Moreover, using convergences \eqref{tlcv}, \eqref{cltrm} and \eqref{tlca}, one can obtain confidence regions for the parameters $v$, $\theta$ and $a$. 
%
%
For instance, for $n=2\,000$, if we denote by $I_{n}(v_{2})$, $I_{n}(\theta_{3})$ and $I_{n}(a_{4})$ the confidence intervals of $v_{2}$, $\theta_{3}$ and $a_{4}$, for a risk $5\%$, one have precisely
\begin{eqnarray*}
I_{n}(v_{2})&=&[-0.1160; 0.4452],\\
I_{n}(\theta_{3})&=&[-0.1211; 0.0226],\\
I_{n}(a_{4})&=&[-2.7225;-2.1224].
\end{eqnarray*}
The lengths of these intervals are respectively 0.5612, 0.1437 and 0.6001. Consequently, the lengths of $I_{n}(v_{2})$, $I_{n}(\theta_{3})$ and $I_{n}(a_{4})$ are small, which confirm the good performance of our parametric estimation procedure.
All these confidence intervals are drawn in Figure \ref{ICparam}. 

\begin{figure}[htb]
\vspace{-2ex}
\begin{minipage}[b] {0.48\linewidth}
\centering 
\centerline {\includegraphics[scale=0.5]{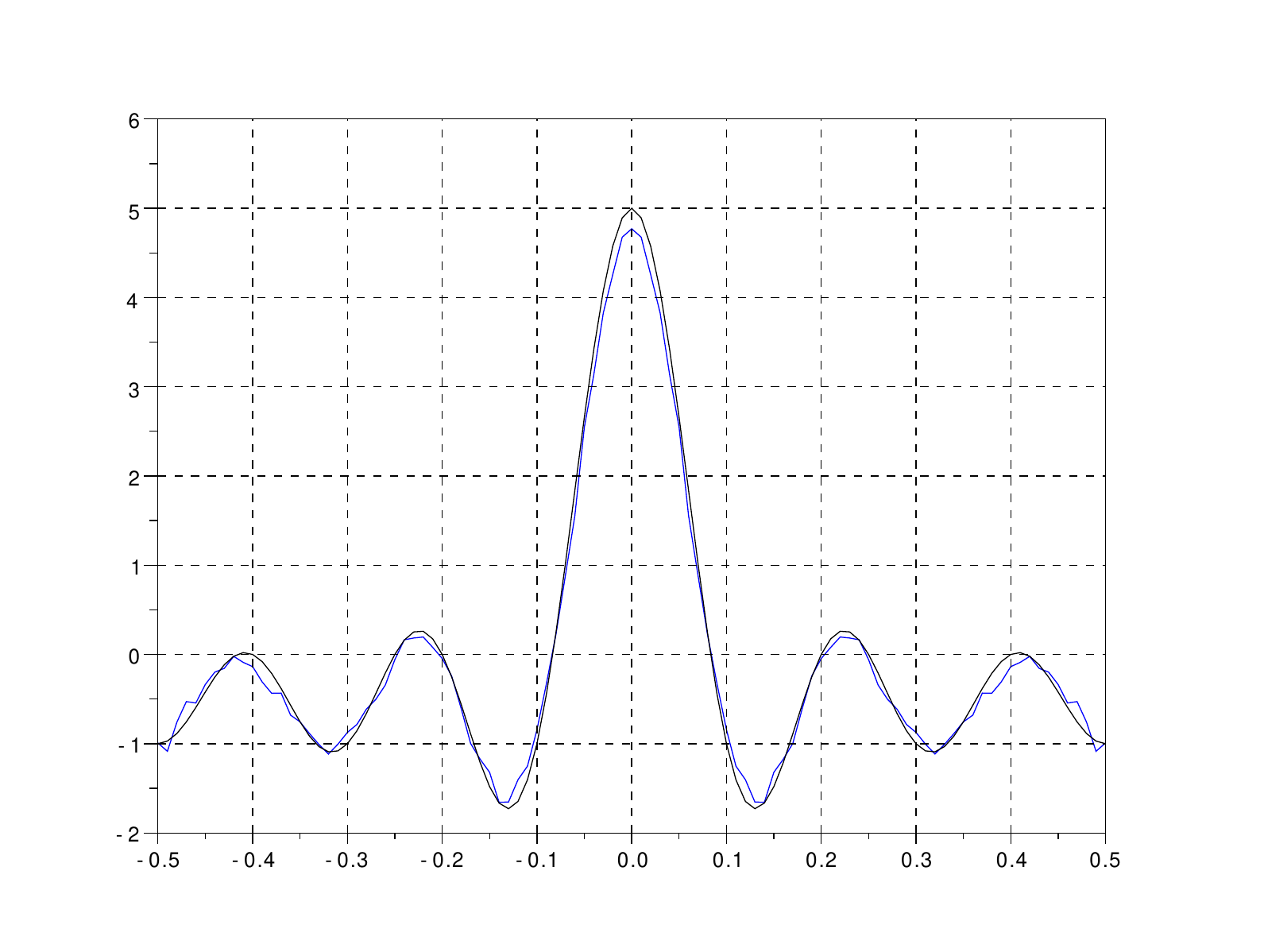} }
\vspace{0.1cm}
\medskip
\end{minipage}
\hfill
\begin{minipage}[b]{0.48\linewidth}
\centering 
\centerline {\includegraphics[scale=0.5]{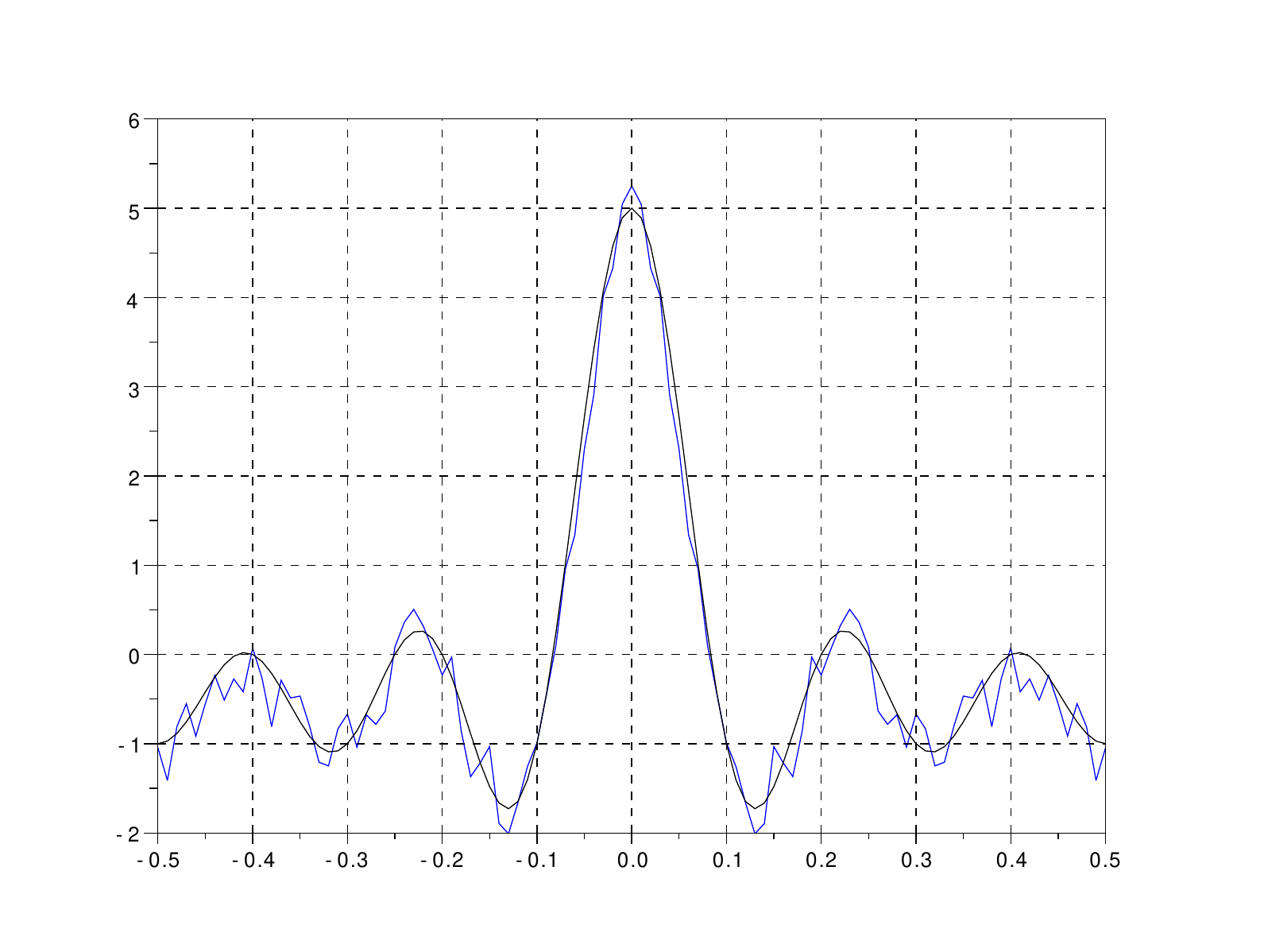} }
\vspace{0.1cm}
\medskip
\end{minipage}
\vspace{-2ex}
\caption{Estimation of $f$ by $\wh{f}_{n}$ and $\wh{f}_{n,1}$}
\label{nonparam}
\end{figure}

\begin{figure}[htb]
\vspace{-2ex}
\begin{minipage}[b] {0.48\linewidth}
\centering 
\centerline {\includegraphics[scale=0.5]{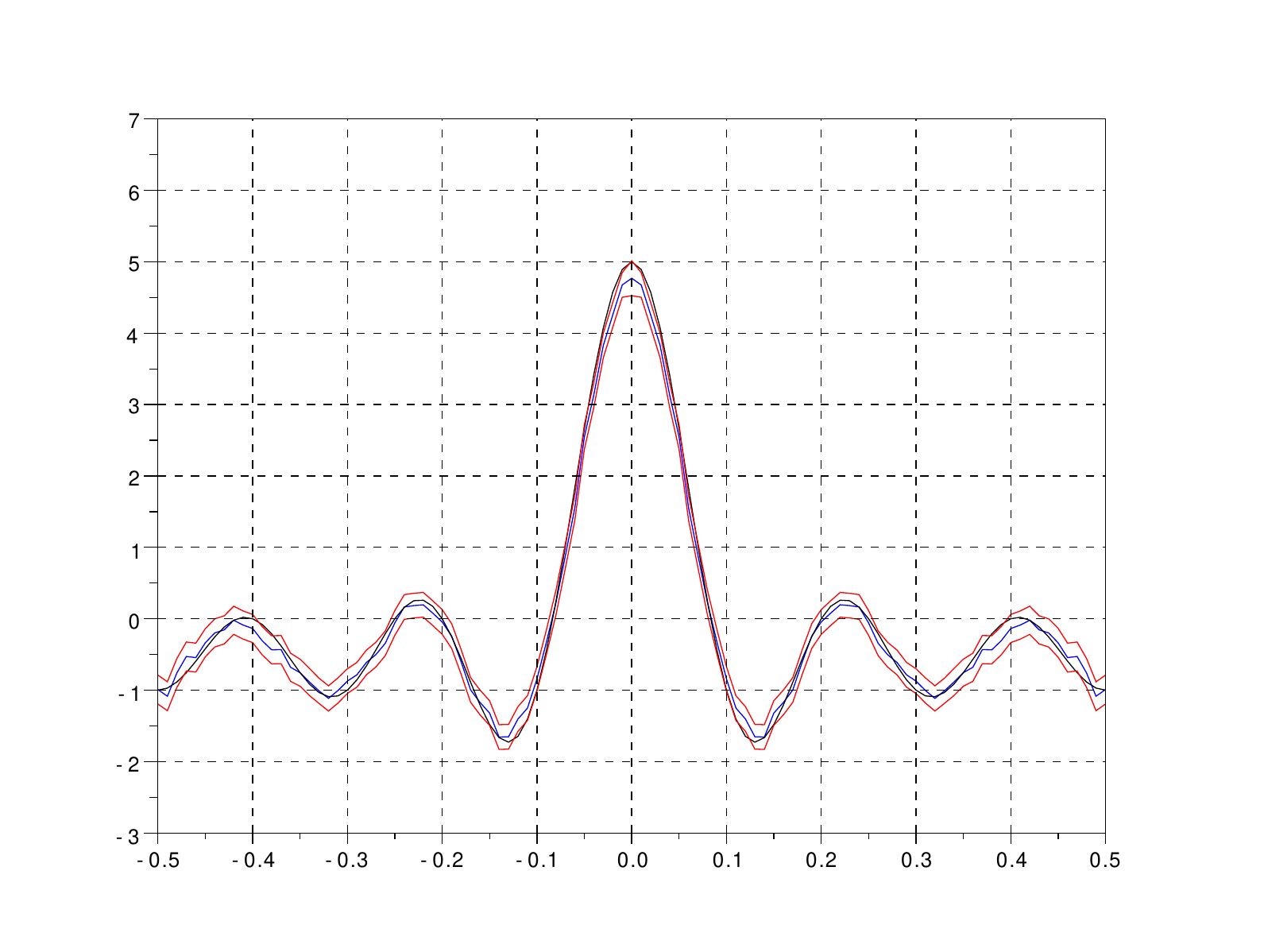} }
\vspace{0.1cm}
\medskip
\end{minipage}
\hfill
\begin{minipage}[b]{0.48\linewidth}
\centering 
\centerline {\includegraphics[scale=0.5]{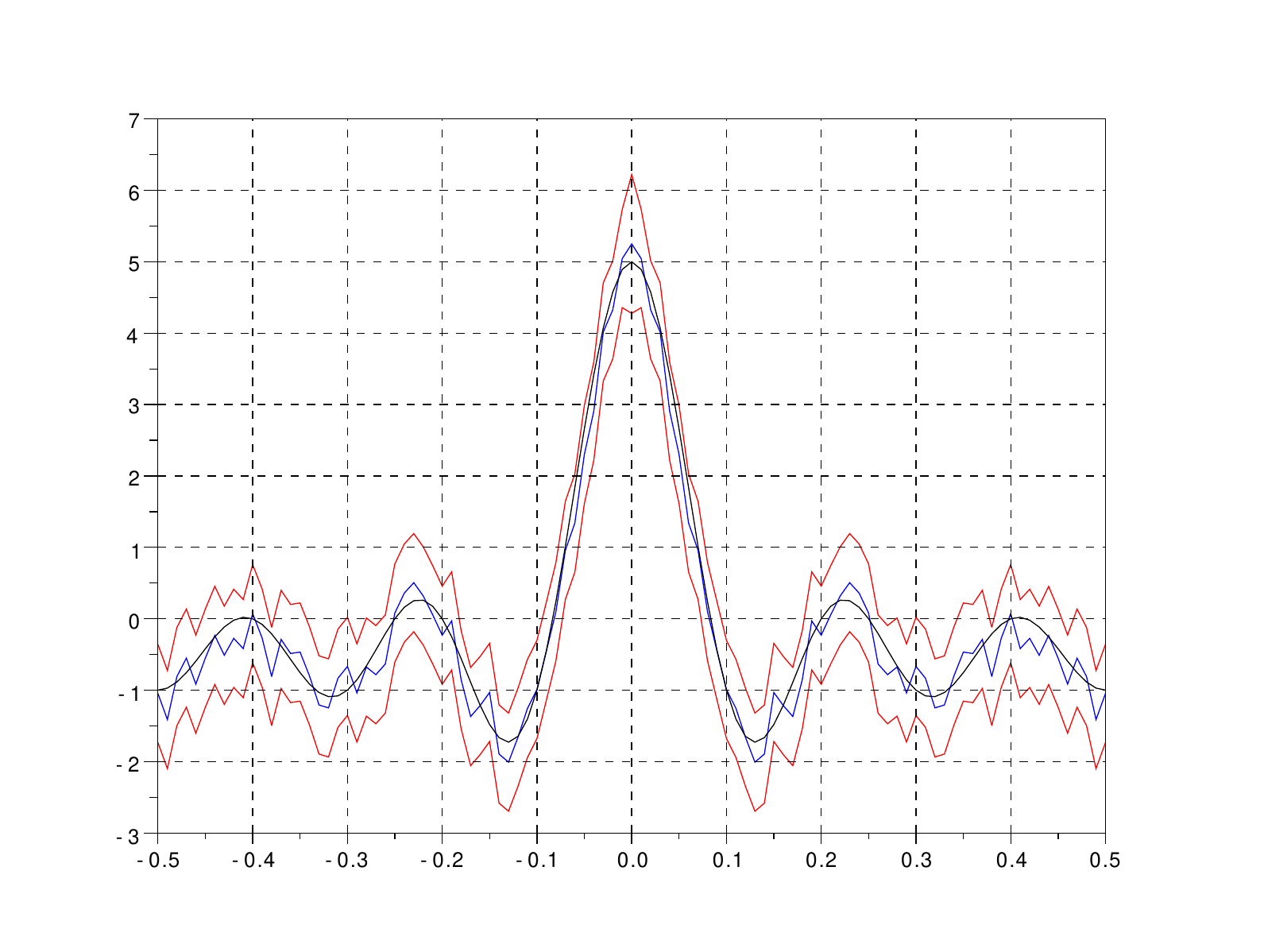} }
\vspace{0.1cm}
\medskip
\end{minipage}
\vspace{-2ex}
\caption{Confidence intervals for $f$}
\label{ICnonparam}
\end{figure}

For the estimation of the regression function $f$, we have chosen $\alpha=9/10$ for the bandwidth sequence $(h_{n})$. Moreover, the kernel $K$ considered is the uniform kernel on $[-1;1]$ and for all $1\leq{j}\leq{p}$,
$\omega_{j}(x)=1/p$. The estimation of the regression function $f$ by $\wh{f}_{n}$ is given on the left side of Figure \ref{nonparam}, and the estimation of $f$ by $\wh{f}_{n,1}$ is given on the right side. 
Furthermore, it follows from convergences \eqref{cltnwx} and \eqref{cltnwzero} that for $n=2\,000$ and for all $x\in{[-1/2;1/2]}$, a confidence interval for $f(x)$ is given by
$$
K_{n}(x)=\left[\wh{f}_{n}(x)-q_{\beta}\frac{\widehat{w}_{n}(x,\wh{\theta}_{n})}{\sqrt{n h_{n}}},\wh{f}_{n}(x)+q_{\beta}\frac{\widehat{w}_{n}(x,\wh{\theta}_{n})}{\sqrt{n h_{n}}}\right]
$$
whereas one can deduce from convergences (3.3) and (3.4) of Theorem 3.2 of \cite{BF10} that for $n=2\,000$ and for all $x\in{[-1/2;1/2]}$, a confidence interval for $f(x)$ is given by
$$
J_{n}(x)=\left[\wh{f}_{n,1}(x)-q_{\beta}\frac{\widehat{v}_{n}(x,\wh{\theta}_{n,1})}{\sqrt{n h_{n}}},\wh{f}_{n,1}(x)+q_{\beta}\frac{\widehat{v}_{n}(x,\wh{\theta}_{n,1})}{\sqrt{n h_{n}}}\right],
$$
where $q_{\beta}$ stands for the quantile of order $0<\beta<1$ of the $\cN(0,1)$ distribution and $\widehat{w}_{n}^{\,2}(x,\wh{\theta}_{n})$ and $\widehat{v}_{n}^{\,2}(x,\wh{\theta}_{n,1})$ are respectively a consistent estimator of the asymptotic variance
$w^{2}(x,\theta)$ in Theorem \ref{thmcltnw} and of $v^{2}(x,\theta_{1})$ in Theorem 3.2 of \cite{BF10}. In our particular case, $\nu^{2}=1/2$, and a numerical calculation leads to
\begin{eqnarray*}
   w^{2}(x,\theta)= \left \{ \begin{array}{lll}
   0.0114 & \ \text{ if } \ -1/2\leq{x}\leq-23/50\hspace{3mm}\text{and}\hspace{3mm}23/50\leq x\leq 1/2, \vspace{1ex} \\
    0.0108& \ \text{ if } \ -23/50<x\leq -9/25\hspace{3mm}\text{and}\hspace{3mm}9/25\leq x<23/50, \vspace{1ex} \\
    0.0099& \ \text{ if } \ -9/25<x\leq -17/50\hspace{3mm}\text{and}\hspace{3mm}17/50\leq x<9/25, \vspace{1ex} \\
    0.0086&\ \text{ if } \  -17/50<x\leq -31/100\hspace{3mm}\text{and}\hspace{3mm}31/100\leq x<17/50, \vspace{1ex} \\
    0.0083&\ \text{ if } \ -31/100<x< 0\hspace{3mm}\text{and}\hspace{3mm}0<x< 31/100, \vspace{1ex} \\
    0.0166&\ \text{ if } x= 0,\vspace{1ex} \\
   \end{array} \nonumber \right.
\end{eqnarray*}
 and
 \begin{eqnarray*}
   v^{2}(x,\theta_{1}) = \left \{ \begin{array}{lll}
    5/38 & \ \text{ if } \ x\neq0, \vspace{1ex} \\
    5/19 & \ \text{ if } \ x=0. \vspace{1ex} \\
   \end{array} \nonumber \right.
\end{eqnarray*}
Roughly speaking, for all $x\in{[-1/2;1/2]}$, the order of the asymptotic variance $v^{2}(x,\theta_{1})$ obtained from the estimate $\wh{f}_{n,1}$ is ten times greater than the order of the asymptotic variance $w^{2}(x,\theta)$ obtained from $\wh{f}_{n}$. In addition, in this case, the optimal variance given in Remark \ref{rempoids} is, for all $x\in{[-1/2;1/2]}$, of order $10^{-3}$ which is anew ten times smaller than $w^{2}(x,\theta)$.
The confidence intervals $K_{n}(x)$ and $J_{n}(x)$ are drawn in red in Figure \ref{ICnonparam}.
One can observe in Figure \ref{nonparam} that the estimate $\wh{f}_{n}$ is closer to the function $f$ than $\wh{f}_{n,1}$. More precisely, the estimation of $f$ by $\wh{f}_{n}$ is better than the one by $\wh{f}_{n,1}$ because the lengths of the confidence intervals $K_{n}(x)$ are smaller than the ones of the confidence intervals $J_{n}(x)$ as one can see in Figure \ref{ICnonparam}, which is due to the order of the two variances $v^{2}(x,\theta_{1})$ and $w^{2}(x,\theta)$. In our case, the largest confidence interval $K_{n}(x)$ is of length $0.3460$ (for $x=0$) and is almost three times smaller than the smallest confidence interval $J_{n}(x)$ which is of length $0.9723$ (for $x=-0.09$ and $x=0.09$). In particular, this justifies the choice of taking a weighed version of the Nadaraya-Watson estimator.

\newpage
\subsection{Modeling of ECG data}
\begin{figure}[htp]
\vspace{-2ex}
\begin{center}
\includegraphics[scale=0.6]{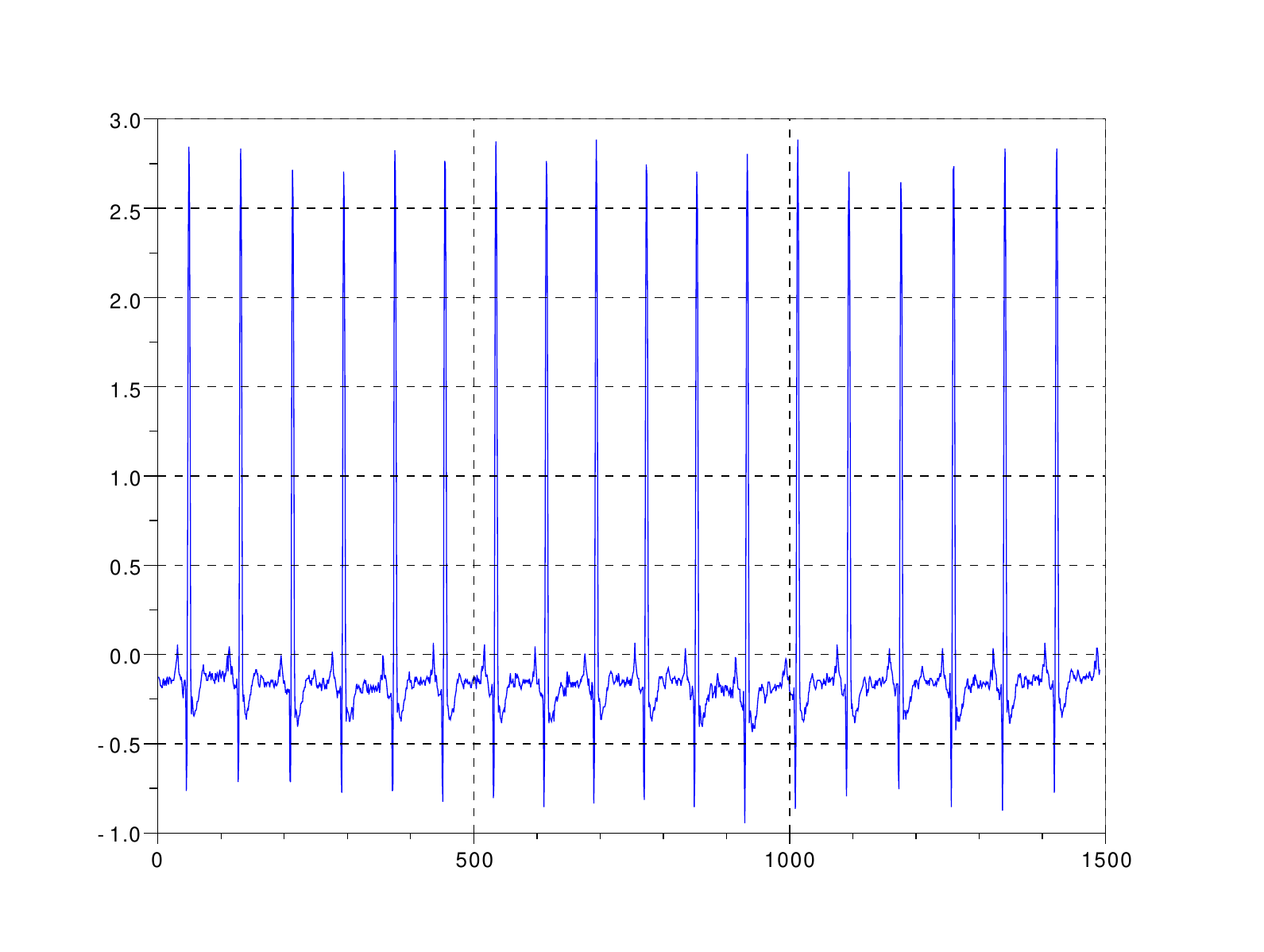}
\vspace{-2ex}
\caption{ECG of a healthy patient}
\label{NormalSignal}
\end{center}
\end{figure}

\begin{figure}[htp]
\vspace{-2ex}
\begin{center}
\includegraphics[scale=0.6]{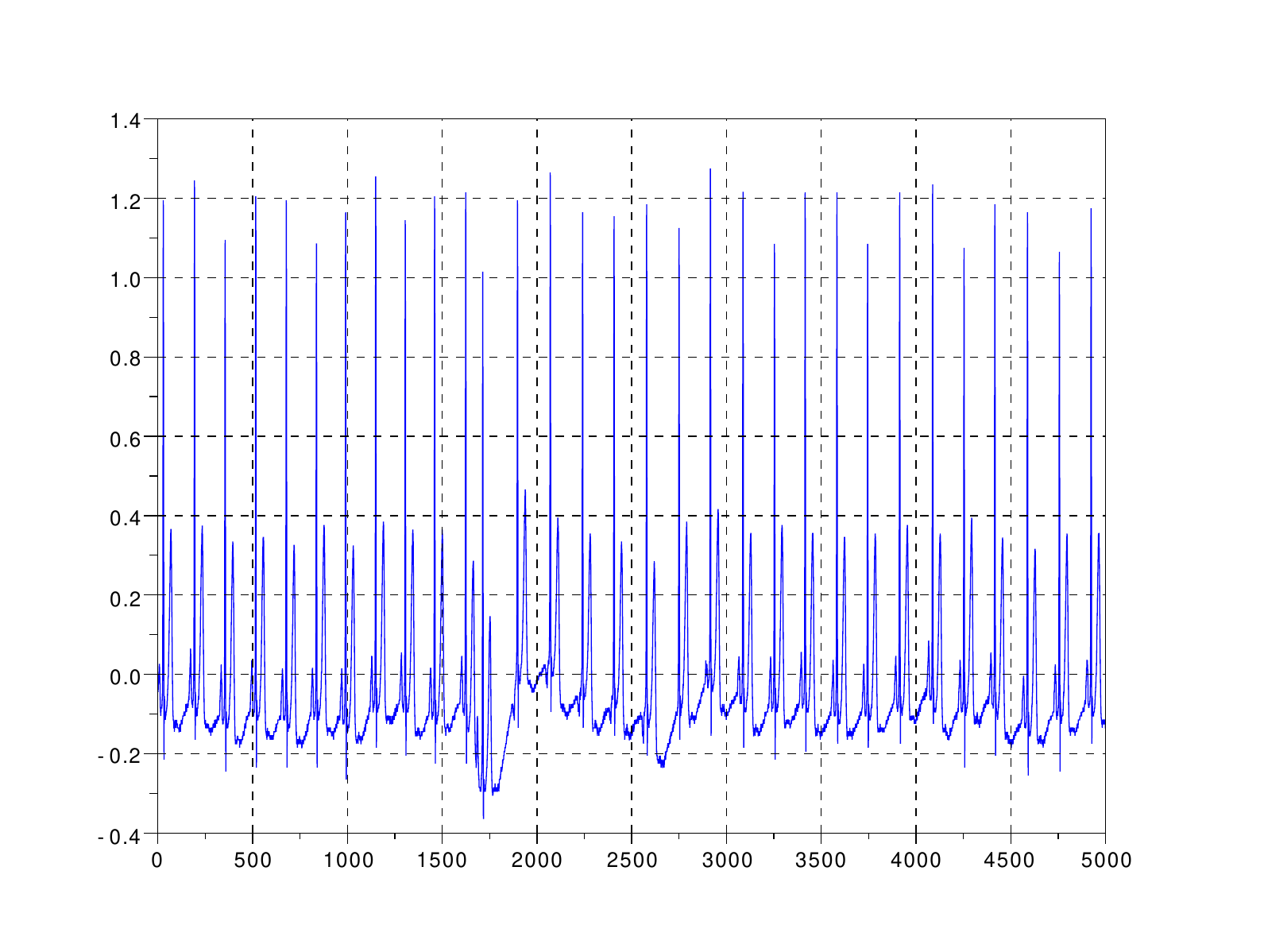}
\vspace{-2ex}
\caption{ECG of a patient having arrythmia}
\label{Signalmalade}
\end{center}
\end{figure}

This section is devoted to real ECG data. An ECG is a recording of the electrical activity of the heart over a period of time, as detected by electrodes attached to the outer surface of the skin.
A typical ECG consists of a P wave, followed by a QRS complex, and a T wave.
Here, we consider two sets of ECG data, one corresponding to a healthy heart and the other to a heart having arrythmia. These data are extracted from the MIT-BIH Database, and they are represented respectively in Figure \ref{NormalSignal} and in Figure \ref{Signalmalade}.

\begin{figure}[htb]
\vspace{-2ex}
\begin{minipage}[b] {0.48\linewidth}
\centering 
\centerline {\includegraphics[scale=0.54]{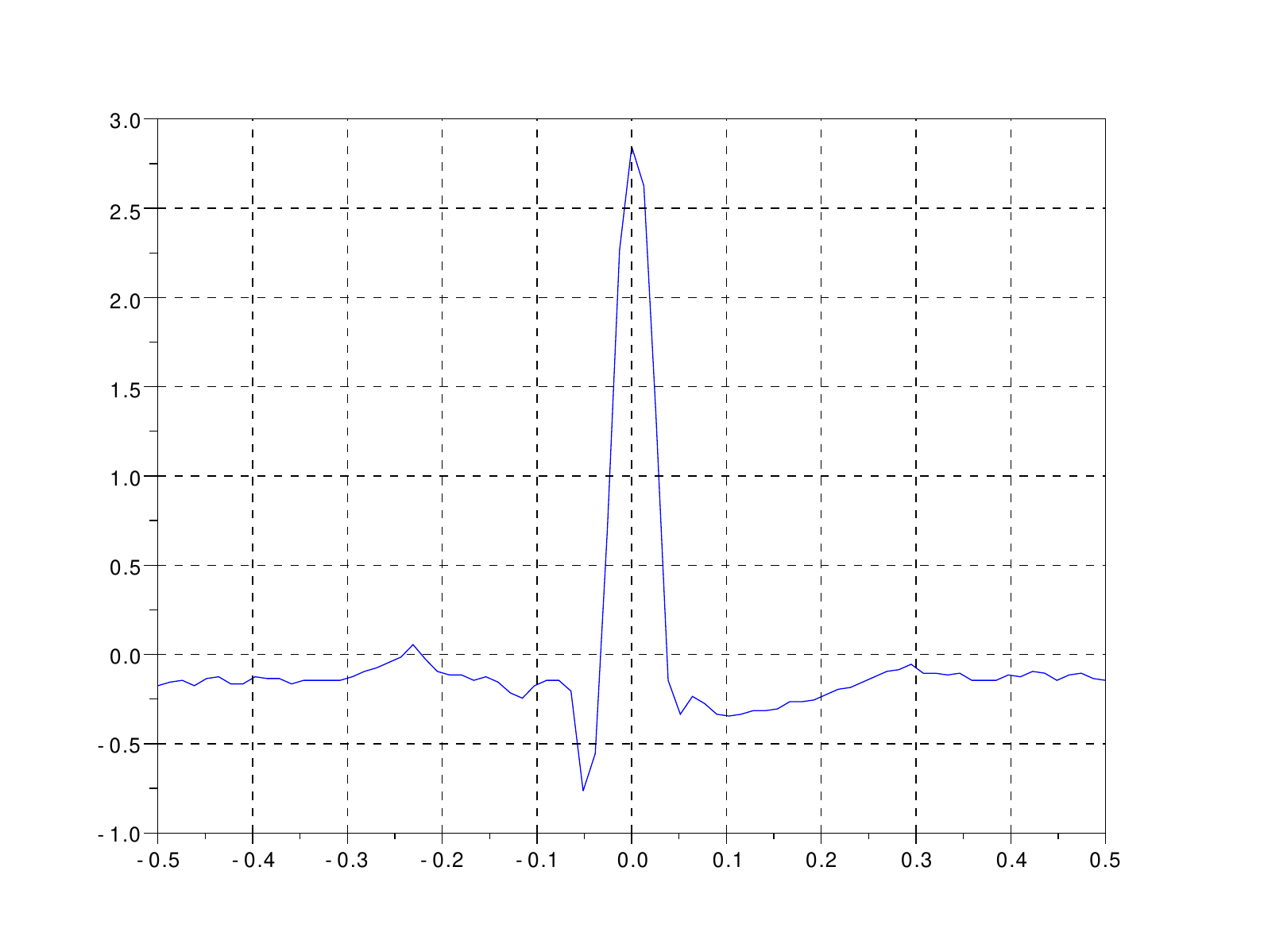} }
\vspace{0.1cm}
\medskip
\end{minipage}
\hfill
\begin{minipage}[b]{0.48\linewidth}
\centering 
\centerline {\includegraphics[scale=0.54]{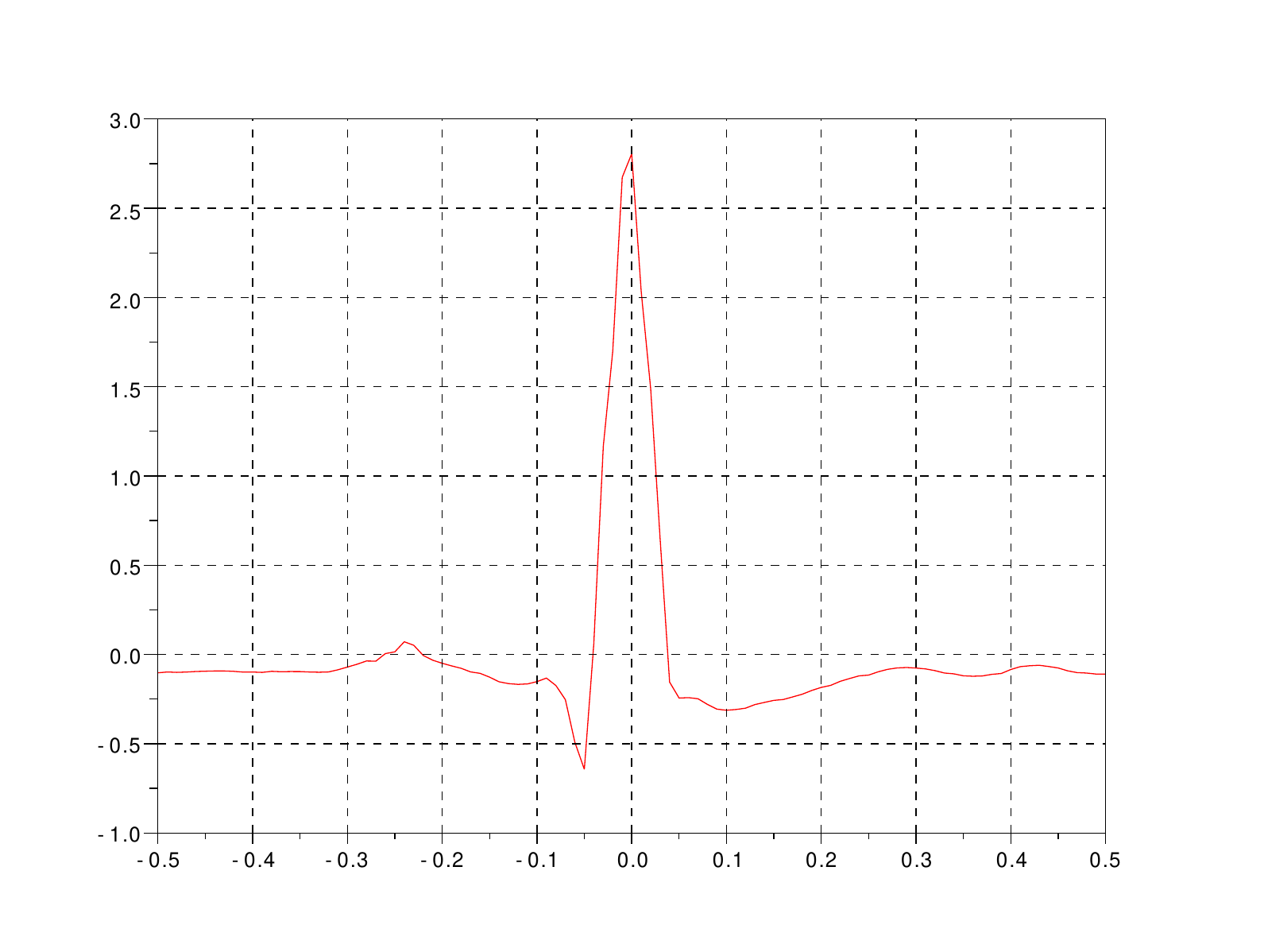} }
\vspace{0.1cm}
\medskip
\end{minipage}
\vspace{-2ex}
\caption{One signal and its reconstruction}
\label{Reconstruction}
\end{figure}

For each ECG recording, one consider that the heart cycle of interest, that is to say the cycle PQRST, is roughly the same at each beat. One also consider that every heart cycle
are noised and that the white noise is independent of the typical shape we want to estimate.
After an appropriate segmentation of the ECGs, one observe signals of same length such that each of them contains an unique PQRST cycle. The segmentation of the ECGs is done by detecting the maximum of each QRS complex and centering the segments around this maxima. It is very important to have segments of same length in order to ensure periodicity. A well-adapted method for the segmentation is the one proposed by Gasser and Kneip \cite{Gasser95}. After the segmentation of the two ECGs, we obtain $p=18$ 
and every segments of length $n=83$ for the healthy heart and $p=15$ and $n=91$ for the ill heart. 
Then, our goal is to estimate the typical shape of each ECG, corresponding to the function $f$ in the model \eqref{Sempar}.
Firstly, the estimation for the healthy heart is going to show that our model is well adapted for the problem of modelling an ECG signal. Indeed, for the healthy heart, a good approximation of the heart cycle is to take the average of the $p$ different signals, whereas the different parameters $v$, $\theta$ and $a$ of the model \eqref{Sempar} are trivial.
The estimation of the typical shape of the ECG for the healthy heart by our procedure is on the right hand-side of Figure \ref{Reconstruction}, whereas the left-side shows one on the original ECG signal. The comparison between the two plots shows that our estimation procedure performs pretty well and that our estimation procedure can be useful for the modeling of an ECG signal.

\begin{figure}[htb]
\vspace{-2ex}
\begin{minipage}[b] {0.48\linewidth}
\centering 
\centerline {\includegraphics[scale=0.54]{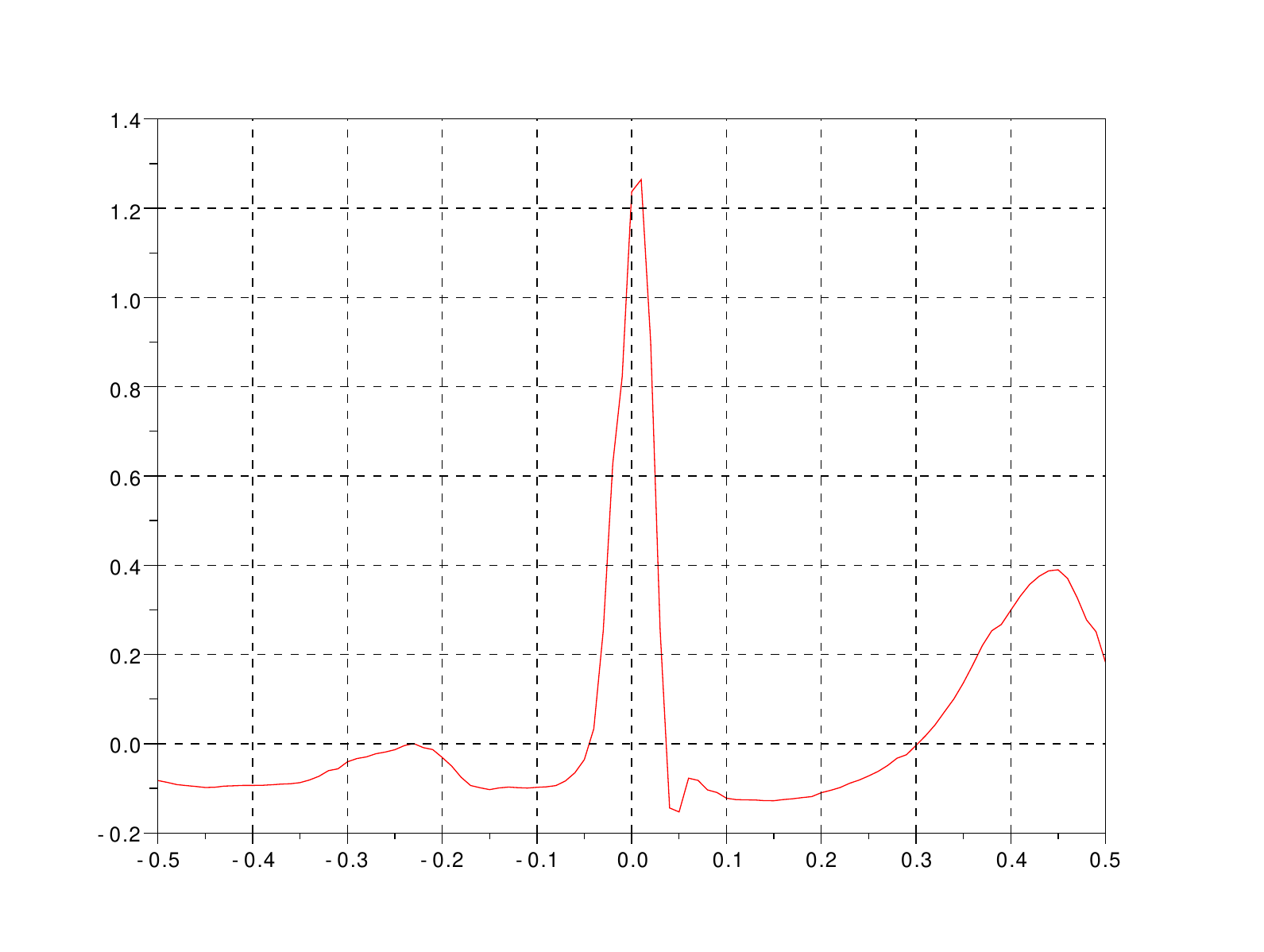} }
\vspace{0.1cm}
\medskip
\end{minipage}
\hfill
\begin{minipage}[b]{0.48\linewidth}
\centering 
\centerline {\includegraphics[scale=0.54]{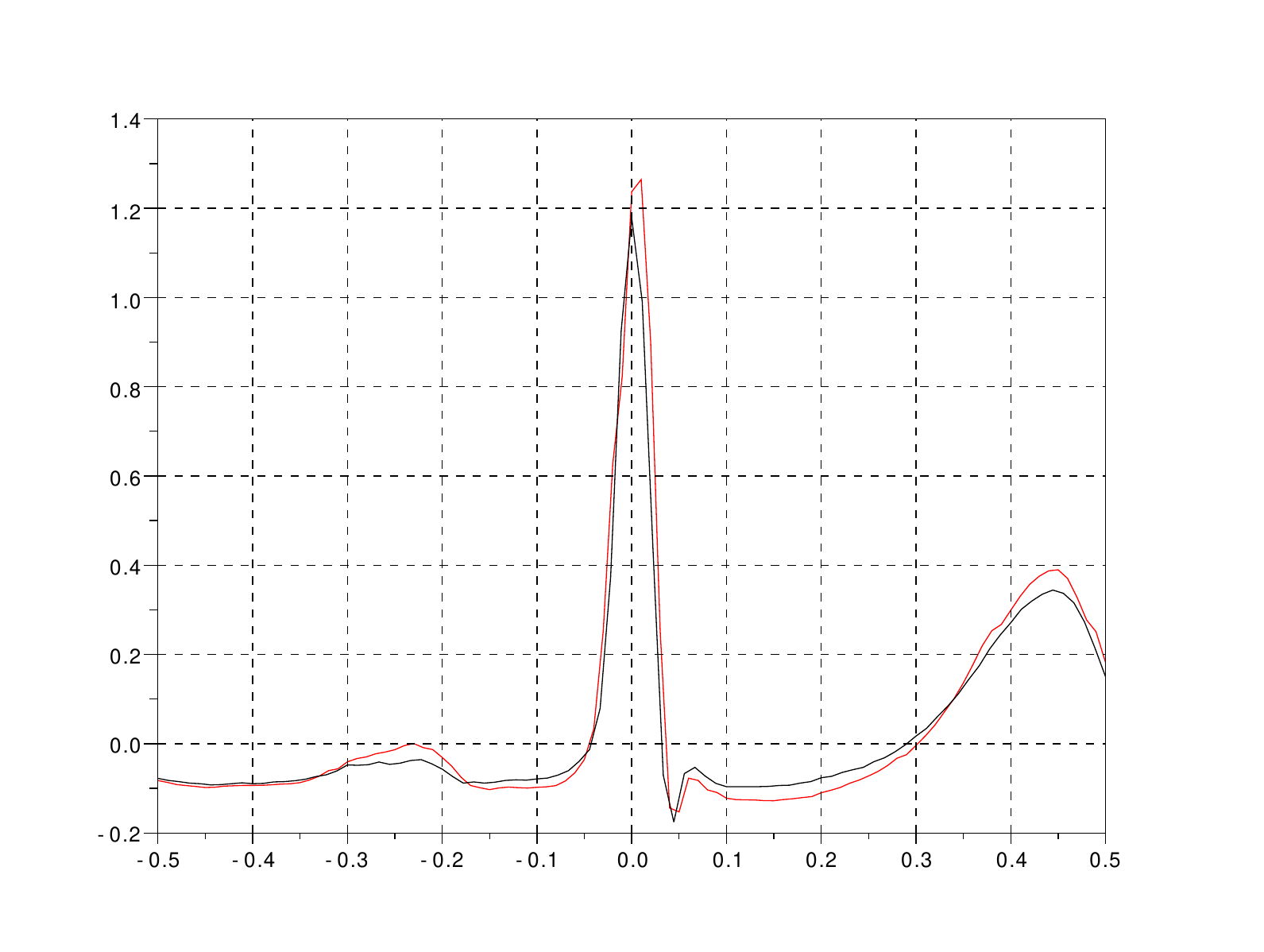} }
\vspace{0.1cm}
\medskip
\end{minipage}
\vspace{-2ex}
\caption{Estimation of $f$ by $\wh{f}_{n}$ (in red) and with the mean average signal (in black)}
\label{Reconstruction2}
\end{figure}

Secondly, we want to estimate the typical shape of the nonhealthy heart, ploted in Figure \ref{Signalmalade}. One see that the electric activity is more irregular than
for the healthy heart, and a simple averaging may lead to a mean cycle that does not correspond to the typical shape of the ECG.
More precisely, we suppose that the model \eqref{Sempar}
$$
Y_{i,j}=a_{j}f(X_{i}-\theta_{j})+v_{j}+\varepsilon_{i,j},
$$
fits the data. The segmentation allows us to have a common shape function $f$ $1$-periodic. Moreover, $f$ is nonsymmetric, but we already saw in Remark \ref{nonsym} that our procedure still works for nonsymetric shape function.
The parameters $a$, $\theta$ and $v$ correspond to the deformation due to the arrythmia relative to the common shape $f$ we want to estimate. For more accuracy, we have to choose one of the $p$ curves as a reference, that is to say where for one $1\leq{j^{*}}\leq{p}$, $a_{j^{*}}=1$, $\theta_{j^{*}}=0$ and $v_{j^{*}}=0$. For this choice, we consider a criterion of residual variance. More precisely, we first consider the model \eqref{Sempar} with $a_{1}=1$, $\theta_{1}=0$ and $v_{1}=0$. From this model, we apply our procedure to estimate the different parameters $a$, $\theta$ and $v$ and the shape function $f$ respectively by $\wh{a}_{n}$, $\wh{\theta}_{n}$, $\wh{v}_{n}$ and $\wh{f}_{n}$. With these estimates, we then calculate the vector 
$\wh{\sigma}_{n,1}^{2}$ whose $j$-th component $(\wh{\sigma}_{n,1}^{2})_{j}$ is defined by
$$
(\wh{\sigma}_{n,1}^{2})_{j}=\frac{1}{n}\sum_{i=1}^{n}\left(Y_{i,j}-\wh{a}_{n,j}\wh{f}_{n}(X_{i}-\wh{\theta}_{n,j})-\wh{v}_{n,j}\right)^{2}.
$$
Then, we make use of the same procedure by changing the curve of reference, and we finally obtain $p$ vectors of length $p$
$$
\wh{\sigma}_{n,1}^{2},\dots,\wh{\sigma}_{n,p}^{2}.
$$
Finally, the choice of the curve of reference for the modeling of the ECG signal is given by taking
$$
j^{*}=\underset{1\leq{j}\leq{p}}\arg\!\min ||\wh{\sigma}_{n,j}^{2}||_{1}
$$
where $||.||_{1}$ corresponds to the $l^{1}$-norm.
Therefore, we model the ECG by
$$
Y_{i,j}=a_{j}f(X_{i}-\theta_{j})+v_{j}+\varepsilon_{i,j}
$$
where $$1\leq{i}\leq{n}\text{, }\hspace{3mm}1\leq{j}\leq{p}\hspace{3mm}\text{ and }\hspace{3mm}a_{j^{*}}=1\text{, }\theta_{j^{*}}=0\text{, }v_{j^{*}}=0.$$
On our data set, the implementation of this method shows that $j^{*}=3$ and $$||\wh{\sigma}_{n,j^{*}}^{2}||_{1}=0.6095.$$
The result for the estimation of the typical shape $f$ is given in Figure \ref{Reconstruction2}. On the right side of Figure \ref{Reconstruction2} one can compare the estimation of $f$ obtained with our estimate $\wh{f}_{n}$ and with the simple average signal. One can observe that our estimation procedure is better because the P wave is well estimated.

\section{PROOFS OF THE PARAMETRIC RESULTS}


\subsection{Proof of Theorem \ref{cvexternalshift}.}
The convergences \eqref{cvpsv} and \eqref{tlcv} follow from the standard law of large numbers and the standard central limit theorem for martingales with independent increments.
Moreover, we immediately deduce \eqref{lfqv} from Theorem 2.1 of \cite{CM2000}.
$\hfill 
\mathbin{\vbox{\hrule\hbox{\vrule height1.5ex \kern.6em
\vrule height1.5ex}\hrule}}$

\subsection{Proof of Theorem \ref{cvpstheta}.}
The result follows from Theorem 2.1 of \cite{BF10}.

$\hfill 
\mathbin{\vbox{\hrule\hbox{\vrule height1.5ex \kern.6em
\vrule height1.5ex}\hrule}}$

%
%

\subsection{Proof of Theorem \ref{thmcltrm}.}
Our goal is to apply Theorem 2.1 page 330 of Kushner and Yin \cite{KushnerYin03}.
First of all, as $\gamma_{n}=1/n$, 
the condition on the decreasing step is satisfied. Moreover, we already saw that $\wh{\theta}_n$ converges almost surely to $\theta$.
Consequently, all the local assumptions of Theorem 2.1 of Kushner and Yin \cite{KushnerYin03} are satisfied. 
In addition, it is not hard to see from \eqref{DefT} that
$$\dE\left[T_{n+1}|\cF_{n}\right]=\phi\left(\wh{\theta}_{n}\right)\hspace{6mm}\textnormal{a.s.}$$
 Moreover, the function $\phi$ is 
continuously differentiable. 
Hence, $\phi(\theta)=0$ and $D\phi(\theta)$ is the square diagonal matrix defined by
 $$D\phi(\theta)=-2\pi f_{1}\mathrm{diag}\left( a_{1},\dots,a_{p}\right).$$
 By noting $I_{p}$ the identity matrix, the condition $4\pi |f_{1}|\underset{1\leq{j}\leq{p}}\min|a_{j}|>1$ implies that
  $$D\phi(\theta)+\frac{1}{2}I_{p}$$ 
is a negative-definite matrix.
Furthermore, we have that, for all $1\leq{k,l}\leq{p}$,
$$
\dE\left[\text{sign}(a_{k}f_{1})T_{n+1,k}\text{sign}(a_{l}f_{1})T_{n+1,l}|\cF_{n}\right]=\varphi(\wh{\theta}_n)_{k,l}\hspace{1cm}\text{a.s.}
$$
which leads to
$$
\lim_{n\rightarrow \infty}\dE\left[\text{sign}(a_{k}f_{1})T_{n+1,k}\text{sign}(a_{l}f_{1})T_{n+1,l}|\cF_{n}\right]=\varphi\left(\theta\right)_{k,l}\hspace{1cm}\text{a.s.}
$$
Consequently, if we are able to prove that the sequence $(W_n)$ given by
$$
W_{n}=\frac{||\wh{\theta}_{n}-\theta||^{2}}{\gamma_{n}}
$$
is tight, then we shall deduce from Theorem 2.1 of  \cite{KushnerYin03} that
$$
\sqrt{n}(\wh{\theta}_{n}-\theta) \liml \cN_{p}(0, \Sigma(\theta))
$$
where for all $1\leq{k,l}\leq{p}$,
$$
\Sigma(\theta)_{k,l}=\varphi(\theta)_{k,l}\int_{0}^{+\infty}\exp\Big(\big(1-2\pi|f_{1}|(|a_{k}|+|a_{l}|)\big)t\Big)\,dt=\frac{\varphi(\theta)_{k,l}}{2\pi|f_{1}|(|a_{k}|+|a_{l}|)-1}.
$$
Therefore, it remains to prove the tightness of the sequence $(W_{n})$. 
Let $(V_{n})$ be the sequence defined, for all $n\geq1$, by
\begin{equation}
\label{defV}
V_{n}=||\wh{\theta}_{n}-\theta||^{2},
\end{equation}
and $T^{\prime}_{n}$ the sequence defined, for all $1\leq{j}\leq{p}$, by
\begin{equation}
\label{defTprime}
T^{\prime}_{n,j}=\textnormal{sign}\left(a_{j}f_{1}\right)T_{n,j}.
\end{equation}
Then, we clearly have
\begin{eqnarray*}
\label{delta}
V_{n+1}&=&||\wh{\theta}_{n+1}-\theta||^{2}\\
&=&||\pi_{K^{p}}\left(\wh{\theta}_{n}+\gamma_{n+1}T^{\prime}_{n+1}\right)-\theta||^{2}\\
&=&||\pi_{K^{p}}\left(\wh{\theta}_{n}+\gamma_{n+1}T^{\prime}_{n+1}\right)-\pi_{K^{p}}\left(\theta\right)||^{2}\\
&\leq&||\wh{\theta}_{n}+\gamma_{n+1}T^{\prime}_{n+1}-\theta||^{2}
\end{eqnarray*}
as $\pi_{K^{p}}=\left(\pi_{K},\dots,\pi_{K}\right)^{T}$ is a Lipschitz function.
It follows that
\begin{eqnarray*}
\label{delta}
V_{n+1}&\leq&V_{n}+\gamma_{n+1}^{2}||T^{\prime}_{n+1}||^{2}+2\gamma_{n+1}<\wh{\theta}_{n}-\theta,T^{\prime}_{n+1}>\hspace{6mm}\text{a.s.}
\end{eqnarray*}
By taking expectation in the previous inequality, we obtain that there exists a constant $M>0$ such that
\begin{equation}
\label{inegV}
\mathbb{E}[V_{n+1}|\mathcal{F}_{n}]\leq V_{n}+\gamma_{n+1}^{2}M+2\gamma_{n+1}<\wh{\theta}_{n}-\theta,\mathbb{E}[T^{\prime}_{n+1}|\mathcal{F}_{n}]>\hspace{6mm}\text{a.s.}
\end{equation}
Moreover, \eqref{DefT} together with \eqref{defTprime} lead to
\begin{equation}
\label{egalTprime}
\mathbb{E}[T^{\prime}_{n+1}|\mathcal{F}_{n}]=S_{p}(a)\phi\left(\wh{\theta}_{n}\right),
\end{equation}
where
$$
S_{p}(a)=\mathrm{diag}\left(\textnormal{sign}\left(a_{1}f_{1}\right),\dots,\textnormal{sign}\left(a_{p}f_{1}\right)\right).
$$
Hence, we deduce from \eqref{inegV} and \eqref{egalTprime} that
\begin{equation}
\label{inegW}
\mathbb{E}[W_{n+1}|\mathcal{F}_{n}]\leq \frac{V_{n}}{\gamma_{n+1}}+\gamma_{n+1}M+2<\wh{\theta}_{n}-\theta,S_{p}(a)\phi\left(\wh{\theta}_{n}\right)>\hspace{6mm}\text{a.s.}
\end{equation}
Moreover, a Taylor expansion of $\phi$ allows us to write
\begin{eqnarray}
\label{Taylorexp}
<\wh{\theta}_{n}-\theta,S_{p}(a)\phi\left(\wh{\theta}_{n}\right)>&=&<\wh{\theta}_{n}-\theta,2\pi f_{1}S_{p}(a)\mathrm{diag}\left(a_{1},\dots,a_{p}\right)\left(\theta-\wh{\theta}_{n}\right)>\nonumber\\
&+&f_{1}<\wh{\theta}_{n}-\theta,S_{p}(a)\mathrm{diag}\left(a_{1},\dots,a_{p}\right)\mathcal{V}\left(\wh{\theta}_{n}\right)\left(\theta-\wh{\theta}_{n}\right)>,\\
\nonumber
\end{eqnarray}
where for all $t\neq{\theta}$,
\begin{equation*}
\mathcal{V}\left(t\right)=\mathrm{diag}\left(\frac{\sin\left(2\pi(\theta_{1}-t_{1})\right)-2\pi(\theta_{1}-t_{1})}{\theta_{1}-t_{1}},\dots,\frac{\sin\left(2\pi(\theta_{p}-t_{p})\right)-2\pi(\theta_{p}-t_{p})}{\theta_{p}-t_{p}}\right).
\end{equation*}
Moreover, the equality
$$
f_{1}S_{p}(a)\mathrm{diag}\left(a_{1},\dots,a_{p}\right)=L(a),
$$
where
$$
L(a)=\mathrm{diag}\left(|f_{1}a_{1}|,\dots,|f_{1}a_{p}|\right),
$$
together with \eqref{Taylorexp} lead to
\begin{eqnarray}
\label{Taylorexp2}
<\wh{\theta}_{n}-\theta,S_{p}(a)\phi\left(\wh{\theta}_{n}\right)>&=&-2\pi\left(\wh{\theta}_{n}-\theta\right)^{T}L(a)\left(\wh{\theta}_{n}-\theta\right)\nonumber\\
&-&\left(\wh{\theta}_{n}-\theta\right)^{T}L(a)\mathcal{V}\left(\wh{\theta}_{n}\right)\left(\wh{\theta}_{n}-\theta\right).\nonumber\\
\nonumber
\end{eqnarray}
Hence, \eqref{inegW} can be rewritten as
\begin{eqnarray}
\label{inegW2}
\mathbb{E}[W_{n+1}|\mathcal{F}_{n}]&\leq& \left(1+\gamma_{n}\right)W_{n}+\gamma_{n}M-4\pi\left(\wh{\theta}_{n}-\theta\right)^{T}L(a)\left(\wh{\theta}_{n}-\theta\right)\\
&&-2\left(\wh{\theta}_{n}-\theta\right)^{T}L(a)\mathcal{V}\left(\wh{\theta}_{n}\right)\left(\wh{\theta}_{n}-\theta\right).
\nonumber
\end{eqnarray}
Moreover, as
$$
L(a)\geq{\underset{1\leq{j}\leq{p}}\min{|a_{j}f_{1}|I_{p}}},
$$
we deduce from \eqref{inegW2} that
\begin{eqnarray}
\label{inegW3}
\mathbb{E}[W_{n+1}|\mathcal{F}_{n}]&\leq&W_{n}+2q\gamma_{n}W_{n}+M\gamma_{n}-2\left(\wh{\theta}_{n}-\theta\right)^{T}L(a)\left(\wh{\theta}_{n}-\theta\right),
\end{eqnarray}
where $$2q=1-4\pi|f_{1}|\underset{1\leq{j}\leq{p}}{\min}|a_{j}|,$$ which means that $q<0$.
By the continuity of the function $\mathcal{V}$, one can find $0<\varepsilon<1/2$ such that, if $||t -\theta||< \varepsilon$, 
\begin{equation}
\label{majv}
\frac{q}{2|f_{1}|\underset{1\leq{j}\leq{p}}{\min}|a_{j}|}I_{p} < \mathcal{V}(t)<0.
\end{equation}
Moreover, let $A_n$ and $B_n$ be the sets $A_{n}=\{||\wh{\theta}_{n}-\theta||\leq \varepsilon\}$
and 
$$
 B_n={\displaystyle \bigcap_{k=m}^n} A_k
$$
with $1\leq m\leq n$. Then, it follows from \eqref{majv} that
\begin{equation}
\label{majvv}
0<-2|f_{1}|\underset{1\leq{j}\leq{p}}{\min}|a_{j}| \mathcal{V}(\wh{\theta}_{n})\rI_{B_{n}}< -\Bigl(\frac{q}{2}\Bigr)I_{p}\rI_{B_{n}}.
\end{equation}
where $I_{p}$ is the identity matrix. Hence, it follows from the conjunction of \eqref{inegW3}
and \eqref{majvv} that for all $n\geq m$,
\begin{eqnarray}
\dE[W_{n+1}\rI_{B_{n}}|\cF_{n}]&\leq & W_{n}\rI_{B_{n}} +  2\gamma_{n}W_{n}q\rI_{B_{n}} -q\gamma_{n}W_{n}\rI_{B_{n}}
+\gamma_{n}M, \nonumber \\
\label{IF}
&\leq & W_{n}\rI_{B_{n}}(1+q\gamma_{n})+\gamma_{n}M.
\end{eqnarray}
Since $B_{n+1}=B_n \cap A_{n+1}$, $B_{n+1}\subset{B_{n}}$, and we obtain by taking the expectation on both sides
of \eqref{IF} that for all $n\geq m$,
\begin{equation}
\label{II}
\dE[W_{n+1}\rI_{B_{n+1}}]\leq (1+q\gamma_{n})\dE[W_{n}\rI_{B_{n}}]+\gamma_{n}M.
\end{equation}
%
Finally, following the same lines as in the proof of Theorem 2.2 in \cite{BF10}, we obtain that for all $\xi>0$, it exists $K>0$ such that for $m$ large enough,
$$
\sup_{n \geq m}\dP(W_{n}>K)<\xi
$$
which implies the tightness of $(W_{n})$ and completes the proof of Theorem \ref{thmcltrm}.
$\hfill 
\mathbin{\vbox{\hrule\hbox{\vrule height1.5ex \kern.6em
\vrule height1.5ex}\hrule}}$

\subsection{Proof of Theorem \ref{thmlilqsl}.}
The vectorial law of iterated logarithm
given by \eqref{lilrm} follows from Theorem 1 of \cite{Pelletier98},
while the vectorial quadratic strong law given by \eqref{qslrm} can be obtained from Theorem 3 of \cite{Pelletier98}.
$\hfill 
\mathbin{\vbox{\hrule\hbox{\vrule height1.5ex \kern.6em
\vrule height1.5ex}\hrule}}$

\subsection{Proof of Theorem \ref{cva}.}
Recall that for all $1\leq{j}\leq{p}$,
$$
\widehat{a}_{n,j}=\frac{1}{nf_{1}}\sum_{i=1}^{n}\frac{\cos(2\pi(X_{i}-\widehat{\theta}_{i-1,j}))}{g(X_{i})}Y_{i,j}
$$
and
$$
\widetilde{a}_{n,j}=\frac{1}{n\wh{f}_{1,n}}\sum_{i=1}^{n}\frac{\cos(2\pi(X_{i}-\widehat{\theta}_{i-1,j}))}{g(X_{i})}Y_{i,j}.
$$
Then, it is clear that
$$
\widehat{a}_{n}=\frac{1}{nf_{1}}\sum_{i=1}^{n}C\left(X_{i},\wh{\theta}_{i-1}\right)Y_{i},
$$
and
$$
\widetilde{a}_{n}=\frac{1}{n\wh{f}_{1,n}}\sum_{i=1}^{n}C\left(X_{i},\wh{\theta}_{i-1}\right)Y_{i},
$$
where
$$
Y_{i}=\begin{pmatrix}
Y_{i,1}\\
\vdots\\
Y_{i,p}
\end{pmatrix}.
$$
We also have the decompositions
\begin{equation}
\label{decompositionANA}
\widehat{a}_{n}-a=\frac{1}{nf_{1}}S_{n}(a)+\frac{1}{nf_{1}}R_{n}(a),
\end{equation}
and
\begin{equation}
\label{decompositionANAtilde}
\widetilde{a}_{n}-a=\frac{1}{n\wh{f}_{1,n}}\left(S_{n}(a)+\left(f_{1}-\wh{f}_{1,n}\right)a\right)+\frac{1}{n\wh{f}_{1,n}}R_{n}(a),
\end{equation}
with
$$
S_{n}(a)=\sum_{i=1}^{n}\left(C(X_{i},\theta)Y_{i}-f_{1}a\right),
$$
and the remainder
$$
R_{n}(a)=\sum_{i=1}^{n}\left(C(X_{i},\wh{\theta}_{i-1})-C\left(X_{i},\theta\right)\right)Y_{i}.
$$
Moreover,
\begin{eqnarray*}
\left(f_{1}-\wh{f}_{1,n}\right)a&=&\sum_{i=1}^{n}\left(f_{1}-\frac{\cos(2\pi X_{i})}{g(X_{i})}Y_{i,1}\right)a,\\
&=&-e_{1}^{T}S_{n}(a)a
\end{eqnarray*}
Hence, we deduce from \eqref{decompositionANAtilde} that
\begin{equation}
\label{decompositionANAtilde2}
\widetilde{a}_{n}-a=\frac{1}{n\wh{f}_{1,n}}M_{p}S_{n}(a)+\frac{1}{n\wh{f}_{1,n}}R_{n}(a),
\end{equation}
where the matrix $M_{p}$ is given by \eqref{defMp}.
In addition, for all $1\leq{j}\leq{p}$,
\begin{equation}
\label{decompoRn}
R_{n,j}(a)=a_{j}R_{n,j}^{1}(a)+v_{j}R_{n,j}^{2}(a)+R_{n,j}^{3}(a),
\end{equation}
where
$$
R_{n,j}^{1}(a)=\sum_{i=1}^{n}\frac{\Delta c_{i,j}}{g(X_{i})}f(X_{i}-\theta_{j}),
\hspace{5mm}
R_{n,j}^{2}(a)=\sum_{i=1}^{n}\frac{\Delta c_{i,j}}{g(X_{i})},
\hspace{5mm}
R_{n,j}^{3}(a)=\sum_{i=1}^{n}\frac{\Delta c_{i,j}}{g(X_{i})}\varepsilon_{i,j},
$$
and
$$
\Delta c_{i,j}=\cos\left(2\pi(X_{i}-\wh{\theta}_{i-1,j})\right)-\cos\left(2\pi(X_{i}-\theta_{j})\right).
$$
Firstly, since
$$
\dE[C(X_{i},\theta)Y_{i}|\mathcal{F}_{i-1}]=f_{1}a,
$$
the sequence $(S_{n}(a))$ is a vectorial martingale with independent increments.
For all $n\geq{1}$, its predictable variation $\langle S(a)\rangle_{n}$ is given by
\begin{eqnarray*}
\langle S(a)\rangle_{n}&=&\sum_{i=1}^{n}\dE[\left(C(X_{i},\theta)Y_{i}-f_{1}a\right)\left(C(X_{i},\theta)Y_{i}-f_{1}a\right)^{T}|\mathcal{F}_{i-1}]\\
&=&\sum_{i=1}^{n}\textnormal{Cov}\left(C(X_{i},\theta)Y_{i}|\mathcal{F}_{i-1}\right)
\end{eqnarray*}
Then, it is clear that
$$
\underset{n\rightarrow{+\infty}}\lim\frac{\langle S(a)\rangle_{n}}{n}=\Gamma(a)\hspace{6mm}\textnormal{a.s.}
$$
where $\Gamma(a)$ is given by \eqref{defgammaa}.\\
\\
Secondly, since
\begin{eqnarray*}
\dE\Big[\frac{\Delta c_{i,j}}{g(X_{i})}|\mathcal{F}_{i-1}\Big]&=&\int_{-1/2}^{1/2}\cos\left(2\pi(x-\wh{\theta}_{i-1,j})\right)dx-\int_{-1/2}^{1/2}\cos\left(2\pi(x-\theta_{j})\right)dx\\
&=&0,
\end{eqnarray*}
the sequence $\left(R_{n,j}^{2}(a)\right)$ is a square integrable martingale whose predictable variation is given by
$$
\langle R^{2}_{j}(a)\rangle_{n}=\sum_{i=1}^{n}\dE\Big[\frac{\Delta c_{i,j}^{2}}{g^{2}(X_{i})}|\mathcal{F}_{i-1}\Big].
$$
Then, as $\cos$ is a Lipschitz function, we have
$$
|\Delta c_{i,j}|\leq{|\wh{\theta}_{i-1,j}-\theta_{j}|}.
$$
Consequently, since $g$ does not vanish on $[-1/2;1/2]$, there exists a constant $C>0$ such that
\begin{equation}
\label{lipschitz}
\dE\Big[\frac{\Delta c_{i,j}^{2}}{g^{2}(X_{i})}|\mathcal{F}_{i-1}\Big]\leq{C\left(\wh{\theta}_{i-1,j}-\theta_{j}\right)^{2}}.
\end{equation}
Then, it follows from \eqref{qslrm} together with the previous inequality \eqref{lipschitz} that
\begin{equation}
\label{variationRn2}
\langle R^{2}_{j}(a)\rangle_{n}=\mathcal{O}\left(\log(n)\right)\hspace{6mm}\textnormal{a.s.}
\end{equation}
Therefore, we deduce from the strong law of large numbers for martingales given e.g. by Theorem 1.3.15 of \cite{Duflo97} that, for all $1\leq{j}\leq{p}$,
\begin{equation}
\label{logR2}
R_{n,j}^{2}(a)=o\left(\log(n)\right)\hspace{6mm}\textnormal{a.s.}
\end{equation}
Moreover, $\left(R_{n,j}^{3}(a)\right)$ is also a square integrable martingale whose predictable variation is given by
$$
\langle R^{3}_{j}(a)\rangle_{n}=\sigma_{j}^{2}\langle R^{2}_{j}(a)\rangle_{n}.
$$
Then, we immediately deduce from \eqref{variationRn2} that
\begin{equation}
\label{variationRn3}
\langle R^{3}_{j}(a)\rangle_{n}=\mathcal{O}\left(\log(n)\right)\hspace{6mm}\textnormal{a.s.},
\end{equation}
and from the strong law of large numbers for martingales that, for all $1\leq{j}\leq{p}$,
\begin{equation}
\label{logR3}
R_{n,j}^{3}(a)=o\left(\log(n)\right)\hspace{6mm}\textnormal{a.s.}
\end{equation}
Afterwards, for all $1\leq{j}\leq{p}$, with the change of variables $u=x-\theta_{j}$,
\begin{eqnarray*}
\dE\Big[\frac{\Delta c_{i,j}}{g(X_{i})}f(X_{i}-\theta_{j})|\mathcal{F}_{i-1}\Big]&=&\int_{-1/2}^{1/2}\left(\cos(2\pi(x-\wh{\theta}_{i-1,j}))-\cos\left(2\pi(x-\theta_{j})\right)\right)f(x-\theta_{j})dx\\
&=&\int_{-1/2-\theta_{j}}^{1/2-\theta_{j}}\left(\cos(2\pi(u+\theta_{j}-\wh{\theta}_{i-1,j}))-\cos\left(2\pi u\right)\right)f(u)du.
\end{eqnarray*}
Then, the elementary trigonometric equality $$\cos(2\pi(u+\theta_{j}-\wh{\theta}_{i-1,j}))=\cos(2\pi u)\cos(2\pi(\theta_{j}-\wh{\theta}_{i-1,j}))-\sin(2\pi u)\sin(2\pi(\theta_{j}-\wh{\theta}_{i-1,j})),$$ the symmetry and the periodicity of the function $f$ lead to
\begin{eqnarray*}
\dE\Big[\frac{\Delta c_{i,j}}{g(X_{i})}f(X_{i}-\theta_{j})|\mathcal{F}_{i-1}\Big]=f_{1}\left(\cos(2\pi(\theta_{j}-\wh{\theta}_{i-1,j}))-1\right)\hspace{6mm}\textnormal{a.s.}
\end{eqnarray*}
Moreover, for all $|x|<1/2$, we have
$$
|\cos(2\pi x)-1|\leq{2\pi^{2}x^{2}},
$$
which implies that
\begin{equation}
\label{majoC}
|\dE\Big[\frac{\Delta c_{i,j}}{g(X_{i})}f(X_{i}-\theta_{j})|\mathcal{F}_{i-1}\Big]|\leq 2\pi^{2}\left(\theta_{j}-\wh{\theta}_{i-1,j}\right)^{2}\hspace{6mm}\textnormal{a.s.}
\end{equation}
Moreover, we have the decomposition
$$
R_{n,j}^{1}(a)=A_{n,j}(a)+B_{n,j}(a),
$$
where
$$
A_{n,j}(a)=\sum_{i=1}^{n}\left(\frac{\Delta c_{i,j}}{g(X_{i})}f(X_{i}-\theta_{j})-\dE\Big[\frac{\Delta c_{i,j}}{g(X_{i})}f(X_{i}-\theta_{j})|\mathcal{F}_{i-1}\Big]\right),
$$
and
$$
B_{n,j}(a)=\sum_{i=1}^{n}\dE\Big[\frac{\Delta c_{i,j}}{g(X_{i})}f(X_{i}-\theta_{j})|\mathcal{F}_{i-1}\Big].
$$
It follows once again from the quadratic strong law \eqref{qslrm} together with \eqref{majoC} that, for all $1\leq{j}\leq{p}$,
\begin{equation}
\label{logB}
B_{n,j}(a)=\mathcal{O}\left(\log(n)\right)\hspace{6mm}\textnormal{a.s.}
\end{equation}
Moreover, for all $1\leq{j}\leq{p}$, $\left(A_{n,j}(a)\right)$ is a square integrable martingale whose predictable variation $\langle A_{j}(a)\rangle_{n}$ satisfy
$$
\langle A_{j}\rangle_{n}\leq{\sum_{i=1}^{n}\dE\Big[\frac{\Delta c_{i,j}^{2}}{g(X_{i})^{2}}f^{2}(X_{i}-\theta_{j})}|\mathcal{F}_{i-1}\Big].
$$
As the shape function $f$ is bounded, we deduce from \eqref{qslrm} together with \eqref{lipschitz} that
$$
\langle A_{j}(a)\rangle_{n}=\mathcal{O}\left(\log(n)\right)\hspace{6mm}\textnormal{a.s.}
$$
Therefore, we can conclude from the strong law of large numbers for martingales that, for all $1\leq{j}\leq{p}$,
\begin{equation}
\label{logA}
A_{n,j}(a)=o\left(\log(n)\right)\hspace{6mm}\textnormal{a.s.}
\end{equation}
Finally, we infer from \eqref{logR2}, \eqref{logR3} together with \eqref{logB} and \eqref{logA} that, for all $1\leq{j}\leq{p}$,
\begin{equation}
\label{logR}
R_{n,j}(a)=\mathcal{O}\left(\log(n)\right)\hspace{6mm}\textnormal{a.s.}
\end{equation}
Hence, one obtain from \eqref{decompositionANA} that
\begin{equation}
\label{decompositionANA2}
\widehat{a}_{n}-a=\frac{1}{nf_{1}}S_{n}(a)+\mathcal{O}\left(\frac{\log(n)}{n}\right)\hspace{6mm}\textnormal{a.s.}
\end{equation}
and from \eqref{decompositionANAtilde2} that
\begin{equation}
\label{decompositionANAtilde3}
\widetilde{a}_{n}-a=\frac{1}{n\wh{f}_{1,n}}M_{p}S_{n}(a)+\mathcal{O}\left(\frac{\log(n)}{n}\right)\hspace{6mm}\textnormal{a.s.}
\end{equation}
Consequently, as $\wh{f}_{1,n}$ converges almost surely to $f_{1}$, \eqref{cvpsa} and \eqref{cvpsatilde} follow from the law of large numbers for martingales and \eqref{tlca} and \eqref{tlcatilde} follow from the central limit theorem for martingales and Slutsky's lemma, while one can obtain \eqref{lfqa} \eqref{lfqatilde} from Theorem 2.1 of \cite{CM2000}.

$\hfill 
\mathbin{\vbox{\hrule\hbox{\vrule height1.5ex \kern.6em
\vrule height1.5ex}\hrule}}$


\section{PROOFS OF THE NONPARAMETRIC RESULTS}


\subsection{Proof of Theorem \ref{thmaspnw}.}

For $x\in{[-1/2;1/2]}$, denote by $\check{f}_{n,j}(x)$ the sequence defined for $n\geq1$ and $1\leq{j}\leq{p}$, by
\begin{eqnarray}
\label{ftilde}
\check{f}_{n,j}(x)=\widehat{a}_{n,j}\wh{f}_{n,j}(x).
\end{eqnarray}
We can rewrite
\begin{eqnarray}
\label{decompftilde}
\check{f}_{n,j}(x)=\check{f}_{n,j}^{1}(x)+\check{f}_{n,j}^{2}(x),
\end{eqnarray}
where
$$
\check{f}_{n,j}^{1}(x)=\frac{\sum_{i=1}^{n}\left(W_{i,j}(x)+W_{i,j}(-x)\right)\left(Y_{i,j}-v_{j}\right)}{\sum_{i=1}^{n}\left(W_{i,j}(x)+W_{i,j}(-x)\right)},
$$
and
$$
\check{f}_{n,j}^{2}(x)=\frac{\sum_{i=1}^{n}\left(W_{i,j}(x)+W_{i,j}(-x)\right)\left(v_{j}-\wh{v}_{i-1,j}\right)}{\sum_{i=1}^{n}\left(W_{i,j}(x)+W_{i,j}(-x)\right)}.
$$
On the one hand, it follows from Theorem 3.1 of \cite{BF10} that for any $x\in{[-1/2;1/2]}$, $$\underset{n\rightarrow +\infty}\lim \check{f}_{n,j}^{1}(x)=a_{j}f(x)\hspace{6mm}\textnormal{a.s.}$$ 
On the other hand, the almost sure convergence of $\wh{v}_{i-1,j}$ to $v_{j}$ as $i$ goes to infinity implies by Toeplitz lemma, that for any $x\in{[-1/2;1/2]}$,
$$
\underset{n\rightarrow +\infty}\lim\check{f}_{n,j}^{2}(x)=0\hspace{6mm}\textnormal{a.s.}
$$
Hence, one can conclude that
\begin{equation} 
\label{cvftilde}
\lim_{n\rightarrow \infty}
\check{f}_{n,j}(x)=a_{j}f(x)\hspace{6mm}\text{a.s.}
\end{equation}
Consequently, as $\widehat{a}_{n,j}$ converges almost surely to $a_{j}\neq{0}$ as $n$ goes to infinity, it follows that
\begin{eqnarray}
\label{cvhatf}
\lim_{n\rightarrow \infty}
\wh{f}_{n,j}(x)=f(x)\hspace{6mm}\text{a.s.}
\end{eqnarray}
Finally, \eqref{RNWS} with \eqref{cvhatf} allow us to conclude the proof of Theorem \ref{thmaspnw}.
$\hfill 
\mathbin{\vbox{\hrule\hbox{\vrule height1.5ex \kern.6em
\vrule height1.5ex}\hrule}}$

\subsection{Proof of Theorem \ref{thmcltnw}.}

We shall now proceed to the proof of the asymptotic normality of $\wh{f}_{n}$.
We have, for all $x\in{[-1/2;1/2]}$,
\begin{eqnarray}
\label{deltafn}
\nonumber
\wh{f}_{n}(x)-f(x)&=&\sum_{j=1}^{p}\omega_{j}(x)\left(\wh{f}_{n,j}(x)-f(x)\right)\\
&=&\sum_{j=1}^{p}\omega_{j}(x)\frac{\cM_{n,j}(x)+\cP_{n,j}(x)+\cQ_{n,j}(x)+\cR_{n,j}(x)+\cS_{n,j}(x)}{n\cG_{n,j}(x)}
\end{eqnarray}
where 
\begin{eqnarray*}
\cG_{n,j}(x)&=&\wh{a}_{n,j}\left(\wh{g}_{n,j}(x)+\wh{g}_{n,j}(-x)\right),\\
\cM_{n,j}(x)&=&M_{n,j}(x)+M_{n,j}(-x), \\
\cP_{n,j}(x)&=&P_{n,j}(x)+P_{n,j}(-x), \\ 
\cQ_{n,j}(x)&=&Q_{n,j}(x)+Q_{n,j}(-x),\\
\cR_{n,j}(x)&=&R_{n,j}(x)+R_{n,j}(-x),\\
\cS_{n,j}(x)&=&S_{n,j}(x)+S_{n,j}(-x),\\
\end{eqnarray*}
with $\wh{g}_{n,j}(x) $, $M_{n,j}(x)$, $P_{n,j}(x)$, $Q_{n,j}(x)$, $R_{n,j}(x)$ and $S_{n,j}(x)$ given by
\begin{eqnarray*}
\wh{g}_{n,j}(x)&=&\frac{1}{n}\sum_{i=1}^{n}W_{i,j}(x), \\
M_{n,j}(x)&=&\sum_{i=1}^{n}W_{i,j}(x)\varepsilon_{i,j}, \\
P_{n,j}(x)&=&\widehat{a}_{n,j}\sum_{i=1}^{n}W_{i,j}(x)\left(f(X_{i}-\wh{\theta}_{i-1,j})-f(x)\right), \\ 
Q_{n,j}(x)&=&\widehat{a}_{n,j}\sum_{i=1}^{n}W_{i,j}(x)\left(f(X_{i}-\theta_{j})-f(X_{i}-\wh{\theta}_{i-1,j})\right),\\
R_{n,j}(x)&=&\left(a_{j}-\widehat{a}_{n,j}\right)\sum_{i=1}^{n}W_{i,j}(x)f(X_{i}-\theta_{j}),\\
S_{n,j}(x)&=&\sum_{i=1}^{n}W_{i,j}(x)\left(v_{j}-\wh{v}_{i-1,j}\right).
\end{eqnarray*}
Firstly, (6.28) of \cite{BF10} together with the almost sure convergence of $\wh{a}_{n}$ to $a$ as $n$ goes to infinity, lead to
\begin{equation}
\label{CvgGn}
\underset{n\rightarrow{+\infty}}\lim\cG_{n,j}(x)=a_{j}\left(g(\theta_{j} +x)+g(\theta_{j} -x)\right)\hspace{6mm}\text{a.s.} 
\end{equation}
In addition, we obtain from (6.32) and (6.35) of \cite{BF10} that, for $\alpha>1/3$,
\begin{eqnarray} 
\label{CvgPn}
\cP^{2}_{n,j}(x)&=&o\left(n^{1+\alpha}\right)\hspace{6mm}\text{a.s.},\\
\label{CvgQn}
\cQ^{2}_{n,j}(x)&=&o\left(n^{1+\alpha}\right)\hspace{6mm}\text{a.s.}\\
\nonumber
\end{eqnarray}
Hence, for all $x\in{[-1/2,1/2]}$, we find that
\begin{equation}
\label{CV1}
\lim_{n\rightarrow \infty}
\sqrt{\frac{h_{n}}{n}}\sum_{j=1}^{p}\omega_{j}(x)\frac{\cP_{n,j}(x)+\cQ_{n,j}(x)}{\cG_{n,j}(x)}=0\hspace{6mm}\text{a.s.}
\end{equation}
Secondly, as the shape function $f$ is bounded, it follows that
\begin{equation*}
R_{n,j}(x)=\mathcal{O}\left(|a_{j}-\widehat{a}_{n,j}|\sum_{i=1}^{n}W_{i,j}(x)\right)\hspace{6mm}\text{a.s.}
\end{equation*}
Hence, 
$$
\sum_{i=1}^{n}W_{i,j}(x)=\mathcal{O}(n)\hspace{6mm}\text{a.s.},
$$
ensures that 
\begin{equation}
\label{grandOR}
R_{n,j}(x)=\mathcal{O}\left(n|a_{j}-\widehat{a}_{n,j}|\right)\hspace{6mm}\text{a.s.}
\end{equation}
In addition, we can deduce from \eqref{decompositionANA2} that, for all $1\leq{j}\leq{p}$,
$$
|a_{j}-\widehat{a}_{n,j}|=\mathcal{O}\left(\sqrt{\frac{\log(n)}{n}}\right)\hspace{6mm}\text{a.s.}
$$
which via \eqref{grandOR}, leads to
\begin{equation}
R_{n,j}(x)=\mathcal{O}\left(\sqrt{n\log(n)}\right)\hspace{6mm}\text{a.s.}
\end{equation}
Consequently,
\begin{equation}
\label{comportR}
\cR_{n,j}^{2}(x)=o(n^{1+\alpha})\hspace{6mm}\text{a.s.}
\end{equation}
Thirdly, we have the following inequality
\begin{equation}
|S_{n,j}(x)|\leq{\Lambda_{n,j}(x)+\Sigma_{n,j}(x)}
\end{equation}
where
$$
\Lambda_{n,j}(x)=\sum_{i=1}^{n}\mathcal{L}_{i,j}\left(W_{i,j}(x)-\dE[W_{i,j}(x)|\mathcal{F}_{i-1}]\right)
$$
and
$$
\Sigma_{n,j}(x)=\sum_{i=1}^{n}\mathcal{L}_{i,j}\dE[W_{i,j}(x)|\mathcal{F}_{i-1}]
$$
with $\mathcal{L}_{i,j}=|v_{j}-\wh{v}_{i-1,j}|$.
We deduce from (6.34) of \cite{BF10} together with the Cauchy-Schwarz inequality and the quadratic strong law given by \eqref{lfqv} that
\begin{equation}
\label{majoSigma}
\Sigma_{n,j}(x)=\mathcal{O}\left(\sqrt{n}\left(\sum_{i=1}^{n}\mathcal{L}_{i,j}^{2}\right)^{1/2}\right)=\mathcal{O}\left(\sqrt{n\log(n)}\right)\hspace{6mm}\text{a.s.}
\end{equation}
Moreover, the sequence $\left(\Lambda_{n,j}(x)\right)$ is a martingale whose predictable variation is given by
$$
\langle\Lambda_{j}(x)\rangle_{n}=\mathcal{O}\left(\sum_{i=1}^{n}\mathcal{L}_{i,j}^{2}\dE[W_{i,j}^{2}(x)|\mathcal{F}_{i-1}]\right)\hspace{6mm}\text{a.s.}
$$
Consequently, we obtain one again from the quadratic strong law \eqref{lfqv} that
\begin{equation}
\label{majoLambda}
\langle\Lambda(x)\rangle_{n,j}=\mathcal{O}\left(n^{\alpha}\sum_{i=1}^{n}\mathcal{L}_{i,j}^{2}\right)=\mathcal{O}\left(n^{\alpha}\log(n)\right)\hspace{6mm}\text{a.s.}
\end{equation}
which allows us to show, from the strong law of large numbers for martingales that for any $\gamma>0$,
\begin{equation}
\label{majoLambda3}
\Lambda_{n,j}^{2}(x)=o\left(n^{\alpha}\log(n)^{2+\gamma}\right)\hspace{6mm}\text{a.s.}
\end{equation}
Therefore,
\begin{eqnarray}
\nonumber
\cS_{n,j}^{2}(x)&=&o\left(n^{\alpha}\log(n)^{2+\gamma}\right)+\mathcal{O}\left(n\log(n)\right)\\
\label{comportS}
&=&o\left(n^{1+\alpha}\right)\hspace{6mm}\text{a.s.}
\end{eqnarray}
We are now in position to study the asymptotic behavior of the dominating term $\cM_{n,j}(x)$.
For all $x\in{[-1/2;1/2]}$ and for all $1\leq{j}\leq{p}$, the sequence $(\cM_{n,j}(x))$ is a square-integrable martingale whose predictable variation is given by
$$
\langle\cM_{j}(x)\rangle_{n}=\sigma_{j}^2\sum_{i=1}^{n}\dE\left[\left(W_{i,j}(x)+W_{i,j}(-x)\right)^2|\cF_{i-1}\right].
$$
We deduce from (6.37) of \cite{BF10} that
we have, for $x\neq0$,
\begin{equation}
\label{cvPQV}
\lim_{n\rightarrow \infty}
\frac{\langle\cM_{j}(x)\rangle_{n}}{n^{1+\alpha}}=\frac{\sigma_{j}^{2}\nu^2}{1+\alpha}\left(g(\theta_{j} +x)+g(\theta_{j} -x)\right).
\end{equation}
and from (6.38) of \cite{BF10} that, for $x=0$,
\begin{equation}
\label{cvPQV0}
\lim_{n\rightarrow \infty}
\frac{\langle\cM_{j}(0)\rangle_{n}}{n^{1+\alpha}}=4\frac{\sigma_{j}^{2}\nu^2}{1+\alpha}g(\theta_j ).
\end{equation}

%
%
\noindent Moreover, according to (6.39) of \cite{BF10}, as $(\varepsilon_{i,j})$ has a moment of order $>2$, Lindeberg condition is satisfied for $\cM_{n,j}(x)$.
We can conclude
from the central limit theorem for martingales given e.g. by Corollary 2.1.10 of \cite{Duflo97}
that for all
$x\in{[-1/2;1/2]}$ with $x \neq 0$,
\begin{equation}
\label{cltMn}
\frac{\cM_{n,j}(x)}{\sqrt{n^{1+\alpha}}}\liml \cN\Bigl(0,\frac{\sigma_{j}^2\nu^2}{1+\alpha}\left(g(\theta_{j} +x)+g(\theta_{j} -x)\right)\Bigr),
\end{equation}
while, for $x=0$,
\begin{equation}
\label{cltMnzero}
\frac{\cM_{n,j}(0)}{\sqrt{n^{1+\alpha}}}\liml \cN\Bigl(0,4\frac{\sigma_{j}^2\nu^2}{1+\alpha}g(\theta_{j})\Bigr).
\end{equation}
Finally, it follows from \eqref{cltMn} and \eqref{cltMnzero} and the independence of $\varepsilon_{i,1},\dots,\varepsilon_{i,p}$ together with the previous convergence \eqref{CvgGn} and Slutsky's theorem that, for all $x\in{[-1/2,1/2]}$ with $x \neq 0$,
\begin{equation}
\label{cltTn}
\frac{1}{\sqrt{n^{1+\alpha}}}\sum_{j=1}^{p}\omega_{j}(x)\frac{\cM_{n,j}(x)}{\cG_{n,j}(x)}\liml \cN\Bigl(0,\frac{ \nu^2}{1+\alpha}\sum_{j=1}^{p}\frac{\sigma_{j}^2\omega_{j}^{2}(x)}{a_{j}^{2}\left(g(\theta_{j} +x)+g(\theta_{j} -x)\right)}\Bigr),
\end{equation}
while, for $x=0$,
\begin{equation}
\label{cltTnzero}
\frac{1}{\sqrt{n^{1+\alpha}}}\sum_{j=1}^{p}\omega_{j}(0)\frac{\cM_{n,j}(0)}{\cG_{n,j}(0)}\liml \cN\Bigl(0,\frac{\nu^2}{1+\alpha}\sum_{j=1}^{p}\frac{\sigma_{j}^2\omega_{j}^{2}(0)}{a_{j}^{2}g(\theta_{j})}\Bigr).
\end{equation}
Then, the conjunction of  \eqref{CV1}, \eqref{comportR}, \eqref{comportS} together with the two previous convergences \eqref{cltTn} and \eqref{cltTnzero} and Slutsky's theorem let us to conclude that, for all $x\in{[-1/2,1/2]}$ with $x \neq 0$,
$$
\sqrt{nh_n}(\wh{f}_{n}(x)-f(x)) \liml \cN\Bigl(0,\frac{ \nu^2}{1+\alpha}\sum_{j=1}^{p}\frac{\sigma_{j}^2\omega_{j}^2(x)}{a_{j}^{2}\left(g(\theta_{j} +x)+g(\theta_{j} -x)\right)}\Bigr),
$$
while, for $x=0$,
$$
\sqrt{nh_n}(\wh{f}_{n}(0)-f(0)) \liml \cN\Bigl(0,\frac{ \nu^2}{1+\alpha}\sum_{j=1}^{p}\frac{\sigma_{j}^2\omega_{j}^2(0)}{a_{j}^{2}g(\theta_{j})}\Bigr),
$$
which completes the proof of Theorem \ref{thmcltnw}.

$\hfill 
\mathbin{\vbox{\hrule\hbox{\vrule height1.5ex \kern.6em
\vrule height1.5ex}\hrule}}$


\section{APPENDIX}

In order to prove every identifiability conditions, let us consider that for a given vector of parameters $\left(a,\theta,v\right)$ satisfying $(\mathcal{H}_4)$ and a given shape function $f$ satisfying $(\mathcal{H}_2)$ and $(\mathcal{H}_3)$, one can find another vector of parameters $\left(a^{*},\theta^{*},v^{*}\right)$ satisfying $(\mathcal{H}_4)$
and an other shape function $f^{*}$ satisfying $(\mathcal{H}_2)$ and $(\mathcal{H}_3)$ such that for all $1\leq{j}\leq{p}$ and for all $x\in{\dR}$, \eqref{identifiability} is true, that is to say
\begin{equation}
\label{equation}
a_{j}f(x-\theta_{j})+v_{j}=a^{*}_{j}f^{*}(x-\theta^{*}_{j})+v^{*}_{j}.
\end{equation}
First of all, the periodicity of $f$ and $(\mathcal{H}_3)$ lead to $v_{j}=v^{*}_{j}$. Then, \eqref{equation} becomes
\begin{equation}
\label{identifiability2}
a_{j}f(x-\theta_{j})=a^{*}_{j}f^{*}(x-\theta^{*}_{j}).
\end{equation}
For $j=1$, the identifiability constraints $a_{1}=a^{*}_{1}$ and $\theta_{1}=\theta^{*}_{1}$ enables us to show that $f(x)=f^{*}(x)$. Then, \eqref{identifiability2} can be rewritten as
\begin{equation}
\label{identifiability3}
a_{j}f(x-\theta_{j})=a^{*}_{j}f(x-\theta^{*}_{j}).
\end{equation}
Hence, by denoting $I_{2}$ the integral of the square of $f$
$$
I_{2}=\int_{0}^{1}f^{2}(x)dx,
$$
and by squaring and integrating into \eqref{identifiability3}, we obtain that
\begin{equation}
\label{identifiability4}
a_{j}^{2}I_{2}=\left(a^{*}_{j}\right)^{2}I_{2},
\end{equation}
leading to 
\begin{equation}
\label{identifiability5}
a_{j}^{2}=\left(a^{*}_{j}\right)^{2}.
\end{equation}
If $a_{j}=a^{*}_{j}$, then $f(x-\theta_{j})=f(x-\theta^{*}_{j})$ and it follows from the constraint $$\underset{1\leq{j}\leq{p}}\max|\theta_{j}|<1/4$$ that $\theta_{j}=\theta_{j}^{*}$. Otherwise, if $a_{j}=-a^{*}_{j}$, then $f(x-\theta_{j})+f(x-\theta^{*}_{j})=0$, which clearly leads to the identitity $f(x)=f(x+2(\theta_{j}-\theta_{j}^{*}))$. Therefore, the constraint $$\underset{1\leq{j}\leq{p}}\max|\theta_{j}|<1/4$$
 implies that $\theta_{j}=\theta_{j}^{*}$. Then, $f(x)=0$, which is impossible. Finally, we have shown that
 $$
 v=v^{*}\text{, }\hspace{3mm}\theta=\theta^{*}\text{, }\hspace{3mm}a=a^{*}\hspace{3mm}\text{ and }\hspace{3mm}f=f^{*},
 $$
 leading to the identifiability of the model \eqref{Sempar}.
 This reasoning fits for the other sets of identifiability constraints.
 $\hfill 
\mathbin{\vbox{\hrule\hbox{\vrule height1.5ex \kern.6em
\vrule height1.5ex}\hrule}}$

\medskip
\noindent \textbf{Acknowledgements.} \textit{The author thanks Bernard Bercu for all his advices and for his thorough readings of the paper.}

\nocite{*}
\bibliographystyle{acm}
\bibliography{RM_SIM_2}

\end{document}